\title{Spectral Measures for $G_2$}
\author{
{\sc David E.\ Evans and Mathew Pugh}\\
 {\footnotesize School of Mathematics, Cardiff University,}\\  {\footnotesize Senghennydd Road, Cardiff CF24 4AG, Wales, U.K.}
}
\date{\today}

\documentclass[12pt]{article}
\usepackage{amsmath,amssymb,amsthm}
\usepackage{graphicx}
\usepackage{framed}
\usepackage{color}
\usepackage[usenames,dvipsnames]{xcolor}
\usepackage{multirow}
\usepackage{cancel}
\usepackage{array}

\newtheorem{Def}{Definition}[section]

\newtheorem{Thm}[Def]{Theorem}

\begin{document}
\maketitle

\begin{abstract}
Spectral measures provide invariants for braided subfactors via fusion modules. In this paper we study joint spectral measures associated to the rank two Lie group $G_2$, including the McKay graphs for the irreducible representations of $G_2$ and its maximal torus, and fusion modules associated to all known $G_2$ modular invariants.
\end{abstract}

{\footnotesize
\tableofcontents
}

\section{Introduction} \label{sect:intro}

Spectral measures associated to the compact Lie group $A_1 = SU(2)$ and its maximal torus, nimrep graphs associated to the $SU(2)$ modular invariants, and the McKay graphs for finite subgroups of $SU(2)$ were studied in \cite{banica/bisch:2007} using information about the generating series of the moments of the spectral measure and the Jones series which is related to the Poincar\'{e} series of the subfactor planar algebra associated to the graph (see e.g. \cite{jones:2001}).
In \cite{evans/pugh:2009v, evans/pugh:2010i} the authors studied spectral measures associated to the compact Lie groups $A_1 = SU(2)$ and $A_2 = SU(3)$ and their maximal tori, nimrep graphs associated to the $SU(2)$ and $SU(3)$ modular invariants, and the McKay graphs for finite subgroups of $SU(2)$ and $SU(3)$, using braided subfactor theory.
Spectral measures associated to the compact rank two Lie groups $B_2$ and $C_2$ are studied in \cite{evans/pugh:2012iv}, and for other compact rank two Lie groups in \cite{evans/pugh:2012iii}.
In this paper and its sequel \cite{evans/pugh:2012ii} we focus on the Lie group $G_2 = \mathrm{Aut}(\mathbb{O})$, the automorphism group of the octonions $\mathbb{O}$. It is a simply connected, compact, real rank two Lie group of dimension 14 and is the smallest of the exceptional Lie groups. It is isomorphic to the subgroup of $SO(7)$ that fixes any particular vector in its 8-dimensional real spinor representation. We determine spectral measures and joint spectral measures for (the adjacency matrices of) various graphs related to the Lie group $G_2$: the McKay (or representation) graphs for the irreducible representations of $G_2$ and its maximal torus $\mathbb{T}^2$, nimrep graphs or fusion modules associated to the $G_2$ modular invariants, and in the sequel \cite{evans/pugh:2012ii}, the McKay graphs for finite subgroups of $G_2$.

Suppose $A$ is a unital $C^{\ast}$-algebra with state $\varphi$.
If $a \in A$ is a self-adjoint operator then there exists a compactly supported probability measure $\nu_a$, the spectral measure of $a$, on the spectrum $\sigma(a) \subset \mathbb{R}$ of $a$, uniquely determined by its moments
$\varphi(a^m) = \int_{\sigma(a)} x^m \mathrm{d}\nu_a (x)$,
for all non-negative integers $m$. Note that $\nu_a$ depends on the choice of state $\varphi$ on $A$.
In the cases we consider, the $C^{\ast}$-algebra $A$ will be a space of operators which act on the following Hilbert space $H$. For the McKay graphs for irreducible representations of $G_2$ we have $H = \ell^2(\mathbb{N}) \otimes \ell^2(\mathbb{N})$ , whilst for the irreducible representations of its maximal torus $\mathbb{T}^2$ we have $H = \ell^2(\mathbb{Z}) \otimes \ell^2(\mathbb{Z})$. For a nimrep graph $\mathcal{G}$ (either associated to a $G_2$ modular invariant or the McKay graph for finite subgroups of $G_2$) with vertex set $\mathcal{G}_0$, we take $H = \ell^2(\mathcal{G}_0)$.

As $G_2$ has rank two, its characters are functions on the maximal torus $\mathbb{T}^2$ of $G_2$. For $SU(2)$ and $SU(3)$ it was convenient to determine the spectral measures for the operator $a$ given by the adjacency matrices of the various graphs related to these groups (the McKay graphs for the irreducible representations of the group and its maximal torus, nimrep graphs associated to modular invariants, and the McKay graphs for finite subgroups) by first determining corresponding measures $\varepsilon_a$ over the maximal tori of $SU(2)$ or $SU(3)$ respectively.
This approach is described as a spectral measure blowup in \cite{banica/bichon:2014}.
In the case of $G_2$, the maximal torus $\mathbb{T}^2$ has dimension one greater than the spectrum $\sigma(a) \subset \mathbb{R}$, so that there is a loss of dimension when passing from the measure $\varepsilon_a$ to $\nu_a$. This means that there is an infinite family of measures $\varepsilon_a$ over $\mathbb{T}^2$ which correspond to the spectral measure $\nu_a$ (details of the relation between $\varepsilon_a$ and $\nu_a$ are given in Section \ref{sect:measures-different_domains}).

In order to remove this ambiguity, we also consider measures over the joint spectrum $\sigma(a,b) \subset \sigma(a) \times \sigma(b) \subset \mathbb{R}^2$ of commuting self-adjoint operators $a$ and $b$. The abelian $C^{\ast}$-algebra $B$ generated by $a$, $b$ and the identity 1 is isomorphic to $C(X)$, where $X$ is the spectrum of $B$. Then the joint spectrum is defined as $\sigma(a,b) = \{ (a(x), b(x)) | \, x \in X \}$. In fact, one can identify the spectrum $X$ with its image $\sigma(a,b)$ in $\mathbb{R}^2$, since the map $x \mapsto (a(x), b(x))$ is continuous and injective, and hence a homeomorphism since $X$ is compact \cite{takesaki:2002}.
In the case where the operators $a$, $b$ act on a finite-dimensional Hilbert space, this is the set of all pairs of real numbers $(\lambda_a,\lambda_b)$ for which there exists a non-zero vector $\phi$ such that $a\phi = \lambda_a \phi$, $b\phi = \lambda_b \phi$.
Then there exists a compactly supported probability measure $\widetilde{\nu}_{a,b}$ on $\sigma(a,b)$, which is uniquely determined by its cross moments
\begin{equation} \label{eqn:cross_moments_sa_operators}
\varphi(a^m b^n) = \int_{\sigma(a,b)} x^m y^n \mathrm{d}\widetilde{\nu}_{a,b} (x,y),
\end{equation}
for all non-negative integers $m$, $n$.
The spectral measure for $a$ is then given by the pushforward $(p_a)_{\ast}(\widetilde{\nu}_{a,b})$ of the joint spectral measure $\widetilde{\nu}_{a,b}$ under the orthogonal projection $p_a$ onto the spectrum $\sigma(a)$.

To study the spectral measures for the nimrep graphs associated to the $G_2$ modular invariants, we use the theory of braided subfactors and $\alpha$-induction, which we now briefly review. For a fuller discussion on braided subfactors and $\alpha$-induction see \cite{bockenhauer/evans:2001, bockenhauer/evans:2002}.
The Verlinde algebra of $G_2$ at level $k$ is represented by a non-degenerately braided system of endomorphisms ${}_N \mathcal{X}_N$ on a type $\mathrm{III}_1$ factor $N$.
Its fusion rules $\{ N_{\lambda \nu}^{\mu} \}$ reproduce exactly those of the positive energy representations of the loop group of $G_2$ at level $k$,
$N_{\lambda} N_{\mu} = \sum_{\nu} N_{\lambda \nu}^{\mu} N_{\nu}$.
Furthermore its statistics generators $S$, $T$ obtained from the braided tensor category ${}_N \mathcal{X}_N$ match exactly those of the Ka\u{c}-Peterson modular $S$, $T$ matrices which perform the conformal character transformations (see footnote 2 in \cite{bockenhauer/evans:2001}).
From the Verlinde formula we see that this family $\{ N_{\lambda} \}$ of commuting normal matrices can be simultaneously diagonalised, i.e.
$N_{\lambda} = \sum_{\sigma} (S_{\sigma, \lambda}/S_{\sigma,0}) S_{\sigma} S_{\sigma}^{\ast}$,
where the summation is over each $\sigma \in {}_N \mathcal{X}_N$ and $0$ is the trivial representation. The intriguing aspect being that the eigenvalues $S_{\sigma, \lambda}/S_{\sigma,0}$ and eigenvectors $S_{\sigma} = \{ S_{\sigma, \mu} \}_{\mu}$ are described by the modular $S$ matrix.

A braided subfactor is an inclusion $N \subset M$ where the dual canonical endomorphism decomposes as a finite combination of elements of the Verlinde algebra, i.e. a finite combination of endomorphisms in ${}_N \mathcal{X}_N$. Such subfactors yield modular invariants through the procedure of $\alpha$-induction which allows two extensions of $\lambda$ on $N$, depending on the use of the braiding or its opposite, to endomorphisms $\alpha^{\pm}_{\lambda} \in {}_M \mathcal{X}_M^{\pm}$ of $M$, so that the matrix $Z_{\lambda,\mu} = \langle \alpha_{\lambda}^+, \alpha_{\mu}^- \rangle$ is a modular invariant \cite{bockenhauer/evans/kawahigashi:1999, bockenhauer/evans:2000, evans:2003}.
The systems ${}_M \mathcal{X}_M^{\pm}$ are called the chiral systems, whilst the intersection ${}_M \mathcal{X}_M^0 = {}_M \mathcal{X}_M^+ \cap {}_M \mathcal{X}_M^-$ is called the neutral system. Then ${}_M \mathcal{X}_M^0 \subset {}_M \mathcal{X}_M^{\pm} \subset {}_M \mathcal{X}_M$, where ${}_M \mathcal{X}_M \subset \mathrm{End}(M)$ denotes a system of endomorphisms consisting of a choice of representative endomorphisms of each irreducible subsector of sectors of the from $[\iota \lambda \overline{\iota}]$, $\lambda \in {}_N \mathcal{X}_N$, where $\iota: N \hookrightarrow M$ is the inclusion map. Although ${}_N \mathcal{X}_N$ is assumed to be braided, the systems ${}_M \mathcal{X}_M^{\pm}$ or ${}_M \mathcal{X}_M$ are not braided in general.
The action of each $N$-$N$ sector $\lambda \in {}_N \mathcal{X}_N$ on the $M$-$N$ sectors ${}_M \mathcal{X}_N$ produces a nimrep (non-negative integer matrix representation of the original Verlinde algebra) $G_{\lambda}$, i.e.
$G_{\lambda} G_{\mu} = \sum_{\nu} N_{\lambda \nu}^{\mu} G_{\nu}$. The spectrum of (each) $G_{\lambda}$ reproduces exactly the diagonal part of the modular invariant \cite{bockenhauer/evans/kawahigashi:2000}.
Since the nimreps are a family of commuting matrices, they can be simultaneously diagonalised and thus the eigenvectors $\psi_{\sigma}$ are the same for the nimrep graphs $G_{\lambda}$ for all $\lambda \in {}_N \mathcal{X}_N$.
We have
$G_{\lambda} = \sum_{\sigma} (S_{\sigma, \lambda}/S_{\sigma,0}) \psi_{\sigma} \psi_{\sigma}^{\ast}$,
where the summation is over each $\sigma \in {}_N \mathcal{X}_N$ with multiplicity given by the modular invariant, i.e. the spectrum of $G_{\lambda}$ is given by $\{ S_{\sigma, \lambda}/S_{\sigma,0}$ with multiplicity $Z_{\sigma,\sigma} \}$. The set of $\mu$ with multiplicity $Z_{\mu,\mu}$ is called the set of exponents of $G$.

Along with the identity invariants for $G_2$ for all levels $k$, there are two exceptional invariants due to conformal embeddings at levels 3, 4 \cite{christe/ravanani:1989} and another exceptional invariant at level 4 \cite{verstegen:1990}. These are all the known $G_2$ modular invariants. Since the centre of $G_2$ is trivial, there are no orbifold modular invariants. This list was shown to be complete for all prime heights $k+4$ such that $k+4 \equiv 5, 7 \, (\mbox{mod } 12)$ \cite{ruelle:1990}, and for all other $k\leq31$ \cite{gannon/ho-kim:1994}.

The paper is organised as follows.
In Section \ref{sect:rep_theoryG2T2} we describe the representation theory of $G_2$ and its maximal torus $\mathbb{T}^2$, and in particular focus on the fundamental representations of $G_2$. In Section \ref{sect:measures-different_domains} we discuss spectral measures for $G_2$ over different domains, showing how measures over a region in the complex plane yields a unique $W$-invariant measure over $\mathbb{T}^2$, where $W$ is the Weyl group of $G_2$.

In Section \ref{sect:measure:A12infty(G2)} we determine the (joint) spectral measures associated to the (adjacency matrices of the) McKay graphs given by the action of the irreducible characters of $G_2$ on its maximal torus $\mathbb{T}^2$, and in Section \ref{sect:measureAinftyG2} the (joint) spectral measures associated to the (adjacency matrices of the) McKay graphs of $G_2$ itself. In both cases we focus on the fundamental representations of $G_2$, and determine these (joint) spectral measures over $\mathbb{T}^2$ and the (joint) spectrum of these adjacency matrices.
Finally in Section \ref{sect:measures_nimrepG2} we determine joint spectral measures over $\mathbb{T}^2$ for nimrep graphs arising from $G_2$ braided subfactors.

\section{Representation theory of $G_2$ and its maximal torus} \label{sect:rep_theoryG2T2}

The irreducible representations $\lambda_{(\mu_1,\mu_2)}$ of $G_2$ are indexed by pairs $(\mu_1,\mu_2) \in \mathbb{N}^2$ such that $\mu_1 \geq \mu_2$.
We denote by $\rho_1 = \lambda_{(1,0)}$ the fundamental representation of $G_2$ of dimension 7, where $\rho_1(G_2) \subset SO(7)$.
The maximal torus of $SO(7)$ is $T = \textrm{diag}(D(\omega_1), D(\omega_2), D(\omega_3), 1)$, for $\omega_i \in \mathbb{T}$, where $D(\omega_i) = \left( \begin{array}{cc} \mathrm{Re}(\omega_i) & -\mathrm{Im}(\omega_i) \\ \mathrm{Im}(\omega_i) & \mathrm{Re}(\omega_i) \end{array} \right)$, and the maximal torus of $G_2$ is the subset of $T$ such that $\omega_1+\omega_2+\omega_3=0$, which is isomorphic to $\mathbb{T}^2$. Then the restriction of $\rho_1$ to $\mathbb{T}^2$ is given by
the $7 \times 7$ block-diagonal matrix
\begin{equation} \label{eqn:restrict_rho1G2_to_T2}
(\rho_1|_{\mathbb{T}^2})(\omega_1,\omega_2) = \textrm{diag}(D(\omega_1), D(\omega_2^{-1}), D(\omega_1^{-1}\omega_2), 1),
\end{equation}
for $(\omega_1,\omega_2) \in \mathbb{T}^2$.
We also denote by $\rho_2 = \lambda_{(1,1)}$ the second fundamental representation of $G_2$, the adjoint representation which has dimension 14.

The Lie group $SU(3)$ is a subgroup of $G_2$. The generating function for the $SU(3) \subset G_2$ branching rules was determined in \cite{gaskell/peccia/sharp:1978}. In particular, the fundamental representations $\rho_1$, $\rho_2$ of $G_2$ branch into the following irreducible $SU(3)$ representations:
\begin{equation} \label{eqn:branch_rulesG2-SU(3)}
\rho_1 \longrightarrow \Sigma_1 \oplus \Sigma_3 \oplus \Sigma_3^{\ast}, \qquad \qquad \rho_2 \longrightarrow \Sigma_3 \oplus \Sigma_3^{\ast} \oplus \Sigma_8.
\end{equation}
The representations $\Sigma_3$, $\Sigma_3^{\ast}$ are conjugate to one another and are the fundamental three-dimensional representations of $SU(3)$. The eight-dimensional representation $\Sigma_8$ is the adjoint representation of $SU(3)$ and is obtained from the product of the fundamental representations by removing one copy of the trivial representation $\Sigma_1$, since $\Sigma_3 \otimes \Sigma_3^{\ast} = \Sigma_1 \oplus \Sigma_8$.
Then from (\ref{eqn:branch_rulesG2-SU(3)}) the restriction of $\rho_2$ to $\mathbb{T}^2$ is given by the $14 \times 14$ block-diagonal matrix
\begin{equation} \label{eqn:restrict_rho2G2_to_T2}
(\rho_2|_{\mathbb{T}^2})(\omega_1,\omega_2)= \textrm{diag}(D(\omega_1), D(\omega_2^{-1}), D(\omega_1^{-1}\omega_2), D(1), D(\omega_1\omega_2), D(\omega_1^2\omega_2^{-1}), D(\omega_1^{-1}\omega_2^2)),
\end{equation}
for $(\omega_1,\omega_2) \in \mathbb{T}^2$.

Let $\{ \chi_{(\mu_1,\mu_2)} \}_{\mu_1,\mu_2 \in \mathbb{N}:\mu_1 \geq \mu_2}$, $\{ \sigma_{(\mu_1,\mu_2)} \}_{\mu_1,\mu_2 \in \mathbb{Z}}$ be the irreducible characters of $G_2$, $\mathbb{T}^2$ respectively, where $\chi_{(\mu_1,\mu_2)} := \chi_{\lambda_{(\mu_1,\mu_2)}}$. The characters $\chi_{(\mu_1,\mu_2)}$ of $G_2$ are self-conjugate and thus are maps from the torus $\mathbb{T}^2$ to an interval $I_{\mu} := \chi_{\mu}(\mathbb{T}^2) \subset \mathbb{R}$.
For $\omega_i \in \mathbb{T}$, $\mu_i \in \mathbb{Z}$, the characters of $\mathbb{T}^2$ are given by $\sigma_{(\mu_1,\mu_2)}(\omega_1,\omega_2) = \omega_1^{\mu_1}\omega_2^{\mu_2}$, and satisfy $\overline{\sigma_{(\mu_1,\mu_2)}} = \sigma_{(-\mu_1,-\mu_2)}$.
If $\sigma_1$, $\sigma_2$ is the restriction of $\chi_1:=\chi_{(1,0)}$, $\chi_2:=\chi_{(1,1)}$ respectively to $\mathbb{T}^2$, we have from (\ref{eqn:restrict_rho1G2_to_T2}) and (\ref{eqn:restrict_rho2G2_to_T2}) that
\begin{align}
\sigma_1 & = \chi_1|_{\mathbb{T}^2} = \sigma_{(0,0)} + \sigma_{(1,0)} + \sigma_{(-1,0)} + \sigma_{(0,-1)} + \sigma_{(0,1)} + \sigma_{(-1,1)} + \sigma_{(1,-1)}, \label{eqn:restriction-chi1} \\
\sigma_2 & = \chi_2|_{\mathbb{T}^2} = \sigma_1 + \sigma_{(0,0)} + \sigma_{(1,1)} + \sigma_{(-1,-1)} + \sigma_{(2,-1)} + \sigma_{(-2,1)} + \sigma_{(1,-2)} + \sigma_{(-1,2)}. \quad \;\; \label{eqn:restriction-chi2}
\end{align}
Then
\begin{equation} \label{eqn:fusion_rules-G2-1}
\sigma_1 \sigma_{(\mu_1,\mu_2)} = \sigma_{(\mu_1,\mu_2)} + \sigma_{(\mu_1+1,\mu_2)} + \sigma_{(\mu_1-1,\mu_2)} + \sigma_{(\mu_1,\mu_2-1)} + \sigma_{(\mu_1,\mu_2+1)} + \sigma_{(\mu_1-1,\mu_2+1)} + \sigma_{(\mu_1+1,\mu_2-1)},
\end{equation}
for any $\mu_1,\mu_2 \in \mathbb{Z}$.

The McKay graph for an irreducible representation $\lambda$ of a group $G$ is the graph whose vertices are labelled by the irreducible representations of $G$, and which, for any two irreducible representations $\lambda_1,\lambda_2$ of $G$, has $N_{\lambda,\lambda_2}^{\lambda_1}$ (directed) edges from $\lambda_1$ to $\lambda_2$, where the decomposition of $\lambda \lambda_1$ into irreducible contains $N_{\lambda,\lambda_2}^{\lambda_1}$ copies of $\lambda_2$.
The McKay graph of $\mathbb{T}^2$ for the first fundamental representation $\rho_1$ is identified with the infinite graph ${}^W \hspace{-2mm} \mathcal{A}^{\rho_1}_{\infty}(G_2)$, illustrated in Figure \ref{fig-A_12infty(G2)1}, whose vertices may be labeled by pairs $(\mu_1,\mu_2) \in \mathbb{Z}^2$ such that there is an edge from $(\mu_1,\mu_2)$ to $(\mu_1,\mu_2)$, $(\mu_1+1,\mu_2)$, $(\mu_1-1,\mu_2)$, $(\mu_1,\mu_2-1)$, $(\mu_1,\mu_2+1)$, $(\mu_1-1,\mu_2+1)$ and $(\mu_1+1,\mu_2-1)$.

\begin{figure}[tb]
\begin{minipage}[t]{7.9cm}
\begin{center}
  \includegraphics[width=35mm]{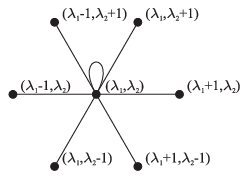}
 \caption{Multiplication by $\rho_1|_{\mathbb{T}^2}$} \label{fig-mult_rho1G2}
\end{center}
\end{minipage}
\hfill
\begin{minipage}[t]{7.9cm}
\begin{center}
  \includegraphics[width=55mm]{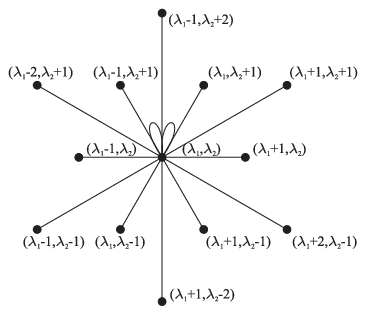}
 \caption{Multiplication by $\rho_2|_{\mathbb{T}^2}$} \label{fig-mult_rho2G2}
\end{center}
\end{minipage}
\end{figure}

\begin{figure}[tb]
\begin{minipage}[t]{7.9cm}
\begin{center}
  \includegraphics[width=75mm]{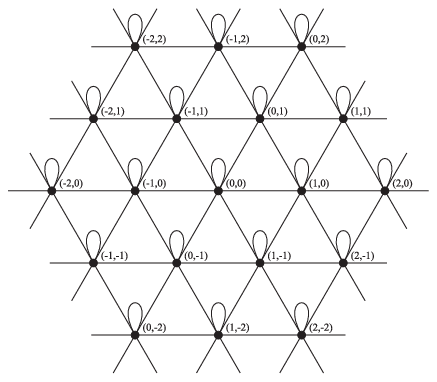}
 \caption{Infinite graph ${}^W \hspace{-2mm} \mathcal{A}^{\rho_1}_{\infty}(G_2)$} \label{fig-A_12infty(G2)1}
\end{center}
\end{minipage}
\hfill
\begin{minipage}[t]{7.9cm}
\begin{center}
  \includegraphics[width=75mm]{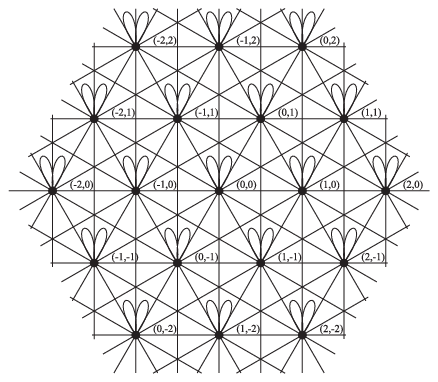}
 \caption{Infinite graph ${}^W \hspace{-2mm} \mathcal{A}^{\rho_2}_{\infty}(G_2)$} \label{fig-A_12infty(G2)2}
\end{center}
\end{minipage}
\end{figure}

Similarly, the McKay graph of $\mathbb{T}^2$ for the second fundamental representation $\rho_2$ is identified with the infinite graph ${}^W \hspace{-2mm} \mathcal{A}^{\rho_2}_{\infty}(G_2)$, illustrated in Figure \ref{fig-A_12infty(G2)2}, where multiplication by $\rho_2$ corresponds to the edges illustrated in Figure \ref{fig-mult_rho2G2}.
These graphs ${}^W \hspace{-2mm} \mathcal{A}^{\rho}_{\infty}(G_2)$ are essentially $W$-unfolded versions of the graphs $\mathcal{A}^{\rho}_{\infty}(G_2)$, where $W$ denotes the Weyl group $D_{12}$ of $G_2$.

By \cite[$\S$3.5]{evans/kawahigashi:1998}, $(\bigotimes_{\mathbb{N}}M_7)^{\mathbb{T}^2} \cong A({}^W \hspace{-2mm} \mathcal{A}^{\rho_1}_{\infty}(G_2))$ and $(\bigotimes_{\mathbb{N}}M_{14})^{\mathbb{T}^2} \cong A({}^W \hspace{-2mm} \mathcal{A}^{\rho_2}_{\infty}(G_2))$.
Here $A(\mathcal{G}) = \overline{\bigcup_k A(\mathcal{G})_k}$ is the path algebra of the graph $\mathcal{G}$, where $A(\mathcal{G})_k$ is the algebra generated by pairs $(\eta_1,\eta_2)$ of paths from the distinguished vertex $\ast$ such that $r(\eta_1)=r(\eta_2)$ and $|\eta_1|=|\eta_2|=k$, with multiplication defined by $(\eta_1,\eta_2) \cdot (\eta_1',\eta_2') = \delta_{\eta_2,\eta_1'}(\eta_1,\eta_2')$.
We now consider instead the fixed point algebra of $\bigotimes_{\mathbb{N}}M_7$, $\bigotimes_{\mathbb{N}}M_{14}$ under the action of the group $G_2$ given by the fundamental representations $\rho_1$, $\rho_2$ respectively, where $G_2$ acts by conjugation on each factor in the infinite tensor product.

The characters $\{ \chi_{(\mu_1,\mu_2)} \}_{\mu_1,\mu_2 \in \mathbb{N}:\mu_1 \geq \mu_2}$ of $G_2$ satisfy
\begin{align*}\chi_1 \chi_{(\mu_1,\mu_2)} & = \left\{ \begin{array}{ll}
\chi_{(\mu_1,\mu_2)} + \chi_{(\mu_1+1,\mu_2)} + \chi_{(\mu_1-1,\mu_2)} + \chi_{(\mu_1,\mu_2-1)} + \chi_{(\mu_1,\mu_2+1)} & \\
\quad + \chi_{(\mu_1-1,\mu_2+1)} + \chi_{(\mu_1+1,\mu_2-1)} & \textrm{ if } \mu_1 \neq \mu_2, \\
\chi_{(\mu_1+1,\mu_2)} + \chi_{(\mu_1,\mu_2-1)} + \chi_{(\mu_1+1,\mu_2-1)} & \textrm{ if } \mu_1 = \mu_2,
\end{array} \right. \\
\chi_2 \chi_{(\mu_1,0)} & = \left\{ \begin{array}{ll}
\chi_1 \chi_{(\mu_1,\mu_2)} + \chi_{(\mu_1+1,1)} + \chi_{(\mu_1-2,1)} + \chi_{(\mu_1-1,2)} & \textrm{ if } \mu_1 \neq \mu_2,\mu_2+1, \\
\chi_{(1,1)} & \textrm{ if } \mu_1 = 0, \\
\chi_{(1,0)} + \chi_{(2,0)} + \chi_{(2,1)} & \textrm{ if } \mu_1 = 1,
\end{array} \right.
\end{align*}
whilst for $\mu_2 \neq 0$,
\begin{align*}
\chi_2 & \chi_{(\mu_1,\mu_2)} \\
& = \left\{ \begin{array}{ll}
\chi_1 \chi_{(\mu_1,\mu_2)} + \chi_{(\mu_1,\mu_2)} + \chi_{(\mu_1-1,\mu_2-1)} + \chi_{(\mu_1+1,\mu_2+1)} & \\
\quad + \chi_{(\mu_1+1,\mu_2-2)} + \chi_{(\mu_1-1,\mu_2+2)} + \chi_{(\mu_1+2,\mu_2-1)} + \chi_{(\mu_1-2,\mu_2+1)} & \textrm{ if } \mu_1 \neq \mu_2,\mu_2+1, \\
\chi_{(\mu_1,\mu_1)} + \chi_{(\mu_1-1,\mu_1-1)} + \chi_{(\mu_1+1,\mu_1+1)} + \chi_{(\mu_1+1,\mu_1-1)} & \\
\quad + \chi_{(\mu_1+1,\mu_1-2)} + \chi_{(\mu_1+2,\mu_1-1)} & \textrm{ if } \mu_1 = \mu_2, \\
2\chi_{(\mu_1,\mu_2)} + \chi_{(\mu_1-1,\mu_2-1)} + \chi_{(\mu_1+1,\mu_2+1)} + \chi_{(\mu_1+1,\mu_2-1)} & \\
\quad + \chi_{(\mu_1+1,\mu_2-2)} + \chi_{(\mu_1+2,\mu_2-1)} + \chi_{(\mu_1+1,\mu_2)} + \chi_{(\mu_1,\mu_2-1)} & \textrm{ if } \mu_1 = \mu_2+1,
\end{array} \right.
\end{align*}
where $\chi_{(\mu_1,\mu_2)} = 0$ if $\mu_2 < 0$ or $\mu_1 < \mu_2$.

\begin{figure}[tb]
\begin{minipage}[t]{7.9cm}
\begin{center}
  \includegraphics[width=70mm]{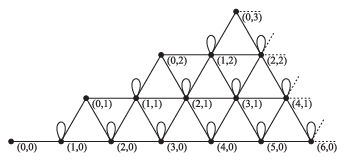}
 \caption{Infinite graph $\mathcal{A}^{\rho_1}_{\infty}(G_2)$} \label{fig-A_infty(G2)1}
\end{center}
\end{minipage}
\hfill
\begin{minipage}[t]{7.9cm}
\begin{center}
  \includegraphics[width=70mm]{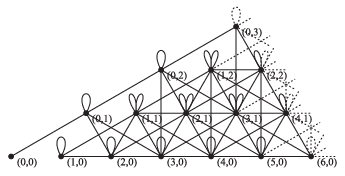}
 \caption{Infinite graph $\mathcal{A}^{\rho_2}_{\infty}(G_2)$} \label{fig-A_infty(G2)2}
\end{center}
\end{minipage}
\end{figure}

Thus the McKay graph of $G_2$ for the first fundamental representation $\rho_1$ is identified with the infinite graph $\mathcal{A}^{\rho_1}_{\infty}(G_2)$, illustrated in Figure \ref{fig-A_infty(G2)1}, where we have made a change of labeling to the Dynkin labels $(\lambda_1,\lambda_2) = (\mu_1-\mu_2,\mu_2)$. This labeling is more convenient in order to be able to define self-adjoint operators $v_N^1$, $v_N^2$ in $\ell^2(\mathbb{N}) \otimes \ell^2(\mathbb{N})$ below. The dashed lines in Figure \ref{fig-A_infty(G2)1} indicate edges that are removed when one restricts to the graph $\mathcal{A}_k(G_2)$ at finite level $k$ (here $k=6$), c.f. Section \ref{sect:measures_AkG2}.

Similarly, the McKay graph of $G_2$ for the second fundamental representation $\rho_2$ is identified with the infinite graph $\mathcal{A}^{\rho_2}_{\infty}(G_2)$, illustrated in Figure \ref{fig-A_infty(G2)2}, again using the Dynkin labels $(\lambda_1,\lambda_2) = (\mu_1-\mu_2,\mu_2)$, and where the dashed lines again indicate edges that are removed when one restricts to the graph $\mathcal{A}_k(G_2)$ at finite level $k$ (here $k=6$).

By \cite[$\S$3.5]{evans/kawahigashi:1998} we have $(\bigotimes_{\mathbb{N}}M_7)^{G_2} \cong A(\mathcal{A}^{\rho_1}_{\infty}(G_2))$ and $(\bigotimes_{\mathbb{N}}M_{14})^{G_2} \cong A(\mathcal{A}^{\rho_2}_{\infty}(G_2))$.

\subsection{Spectral measures over different domains} \label{sect:measures-different_domains}

The Weyl group $W$ of $G_2$ is the dihedral group $D_{12}$ of order 12.
If we consider $D_{12}$ as the subgroup of $GL(2,\mathbb{Z})$ generated by the matrices $T_2$, $T_6$, of orders 2, 6 respectively, given by
\begin{equation} \label{T2,T6}
T_2 = \left( \begin{array}{cc} 0 & -1 \\ -1 & 0 \end{array} \right), \qquad T_6 = \left( \begin{array}{cc} 0 & 1 \\ -1 & 1 \end{array} \right),
\end{equation}
then the action of $D_{12}$ on $\mathbb{T}^2$ given by $T(\omega_1,\omega_2) = (\omega_1^{a_{11}}\omega_2^{a_{12}},\omega_1^{a_{21}}\omega_2^{a_{22}})$, for $T = (a_{il}) \in D_{12}$, leaves $\chi_{\mu}(\omega_1,\omega_2)$ invariant, for any $\mu \in P_+$.
Any $D_{12}$-invariant measure $\varepsilon_{\mu}$ on $\mathbb{T}^2$ yields a pushforward probability measure $\nu_{\mu} = (\chi_{\mu})_{\ast}(\varepsilon_{\mu})$ on $I_{\mu} = \chi_{\mu}( \mathbb{T}^2)\subset \mathbb{R}$ by
\begin{equation} \label{eqn:measures-T2-Ij_G2}
\int_{I_{\mu}} \psi(x) \mathrm{d}\nu_{\mu}(x) = \int_{\mathbb{T}^2} \psi(\chi_{\mu}(\omega_1,\omega_2)) \mathrm{d}\varepsilon_{\mu}(\omega_1,\omega_2),
\end{equation}
for any continuous function $\psi:I_{\mu} \rightarrow \mathbb{C}$, where $\mathrm{d}\varepsilon_{\mu}(\omega_1,\omega_2) = \mathrm{d}\varepsilon_{\mu}(g(\omega_1,\omega_2))$ for all $g \in D_{12}$.
There is a loss of dimension here, in the sense that the integral on the right hand side is over the two-dimensional torus $\mathbb{T}^2$, whereas on the right hand side it is over the interval $I_{\mu}$.
Thus the preimage $I_{\mu}^{-1}[x]$ of any point $x$ in the interior of $I_{\mu}$ is infinite, and there is an infinite family of pullback measures $\varepsilon_{\mu}$ over $\mathbb{T}^2$ for any measure $\nu_{\mu}$ on $I_{\mu}$, that is, any $\varepsilon_{\mu}$ such that $\varepsilon_{\mu}(I_{\mu}^{-1}[x]) = \nu_{\mu}(x)$ for all $x \in I_{\mu}$ will yield the probability measure $\nu_{\mu}$ on $I_{\mu}$ as a pushforward measure by (\ref{eqn:measures-T2-Ij_G2}).
We introduce below an intermediate probability measure $\widetilde{\nu}_{\lambda,\mu}$ which lives over a (two-dimensional) subregion $\mathfrak{D}_{\lambda,\mu} \subset I_{\lambda} \times I_{\mu} \subset \mathbb{R}^2$, for $\lambda,\mu \in P_{++}$, for which there is a unique $D_{12}$-invariant measure $\varepsilon_{\lambda,\mu}$ on $\mathbb{T}^2$. This measure $\widetilde{\nu}_{\lambda,\mu}$ specializes to the spectral measures $\nu_{\lambda}$, $\nu_{\mu}$ of $\lambda$, $\mu$ respectively.

\begin{figure}[tb]
\begin{minipage}[t]{7.9cm}
\begin{center}
  \includegraphics[width=55mm]{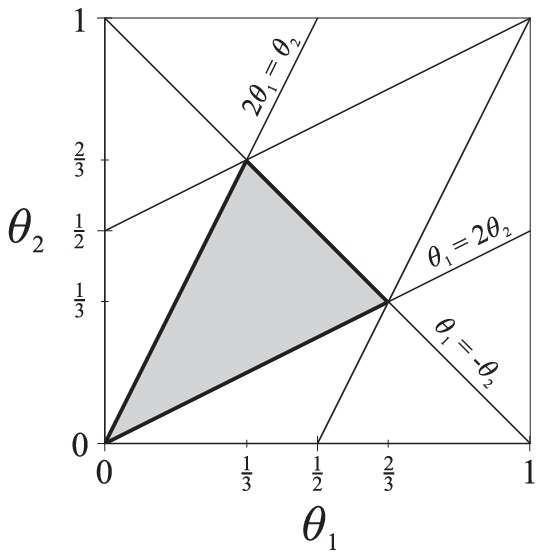}\\
 \caption{A fundamental domain of $\mathbb{T}^2/S_3$.} \label{fig:fund_domain-A2inT2}
\end{center}
\end{minipage}
\hfill
\begin{minipage}[t]{6.9cm}
\begin{center}
  \includegraphics[width=55mm]{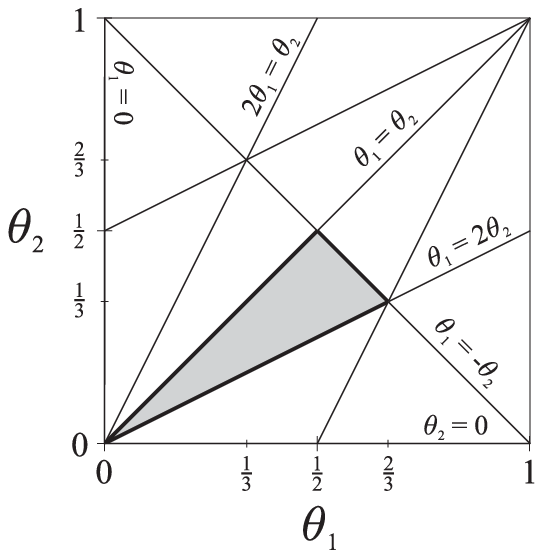}\\
 \caption{A fundamental domain $F$ \mbox{of $\mathbb{T}^2/D_{12}$.}} \label{fig:fund_domain-G2inT2}
\end{center}
\end{minipage}
\end{figure}

The permutation group $S_3$ appears as the subgroup generated by $T_2$ and the matrix $T_6^4 = -T_6$ of order 3 (c.f. \cite[equation (37)]{evans/pugh:2009v}). Then $D_{12}$ is generated by $S_3$ and the $2 \times 2$ matrix $-I$ which sends $\theta_l \leftrightarrow -\theta_l$, $l=1,2$. A fundamental domain of $\mathbb{T}^2/D_{12}$ is thus given by a quotient of the fundamental domain of $\mathbb{T}^2/S_3$, illustrated in Figure \ref{fig:fund_domain-A2inT2} (see \cite{evans/pugh:2009v}), by the $\mathbb{Z}_2$-action given by the matrix $-I$. A fundamental domain $F$ of $\mathbb{T}^2$ under the action of the dihedral group $D_{12}$ is illustrated in Figure \ref{fig:fund_domain-G2inT2}, where the axes are labelled by the parameters $\theta_1$, $\theta_2$ in $(e^{2 \pi i \theta_1},e^{2 \pi i \theta_2}) \in \mathbb{T}^2$. In Figure \ref{fig:fund_domain-G2inT2}, the lines $\theta_1=0$ and $\theta_2=0$ are also boundaries of copies of the fundamental domain $F$ under the action of $D_{12}$, whereas they are not boundaries of copies of the fundamental domain under the action of $S_3$ in Figure \ref{fig:fund_domain-A2inT2}. The torus $\mathbb{T}^2$ contains 12 copies of $F$, so that
\begin{equation} \label{eqn:measureT2=12C}
\int_{\mathbb{T}^2} \phi(\omega_1,\omega_2) \mathrm{d}\varepsilon_{\lambda,\mu}(\omega_1,\omega_2) = 12 \int_{F} \phi(\omega_1,\omega_2) \mathrm{d}\varepsilon_{\lambda,\mu}(\omega_1,\omega_2),
\end{equation}
for any $D_{12}$-invariant function $\phi:\mathbb{T}^2 \rightarrow \mathbb{C}$. The fixed points of $\mathbb{T}^2$ under the action of $S_3$ are the points $(1,1)$, $(e^{2\pi i/3},e^{4\pi i/3})$ and $(e^{4\pi i/3},e^{2\pi i/3})$, but only the point $(1,1)$ is fixed under the action of the whole of $D_{12}$. Under $\chi_{\rho_i}$ the point $(1,1)$ maps to 7, 14 in the intervals $I_{\rho_1}$, $I_{\rho_2}$ respectively, whilst the points $(e^{2\pi i/3},e^{4\pi i/3})$, $(e^{4\pi i/3},e^{2\pi i/3})$ both map to -2 (in both intervals).

Let $x_{\lambda} = \chi_{\lambda}(\omega_1,\omega_2)$ and let $\Psi_{\lambda,\mu}$ be the map $(\omega_1,\omega_2) \mapsto (x_{\lambda},x_{\mu})$. We denote by $\mathfrak{D}_{\lambda,\mu}$ the image of $\Psi_{\lambda,\mu}(F) \; (= \Psi_{\lambda,\mu}(\mathbb{T}^2))$ in $\mathbb{R}^2$.
The joint spectral measure $\widetilde{\nu}_{\lambda,\mu}$ is the measure on $\mathfrak{D}_{\lambda,\mu}$ uniquely determined by its cross-moments as in (\ref{eqn:cross_moments_sa_operators}).
Then there is a unique $D_{12}$-invariant pullback measure $\varepsilon_{\lambda,\mu}$ on $\mathbb{T}^2$ such that
$$\int_{\mathfrak{D}_{\lambda,\mu}} \psi(x_{\lambda},x_{\mu}) \mathrm{d}\widetilde{\nu}_{\lambda,\mu}(x_{\lambda},x_{\mu}) = \int_{\mathbb{T}^2} \psi(\chi_{\lambda}(\omega_1,\omega_2),\chi_{\mu}(\omega_1,\omega_2)) \mathrm{d}\varepsilon_{\lambda,\mu}(\omega_1,\omega_2),$$
for any continuous function $\psi:\mathfrak{D}_{\lambda,\mu} \rightarrow \mathbb{C}$.

Any probability measure on $\mathfrak{D}_{\lambda,\mu}$ yields a probability measure on the interval $I_{\lambda}$, given by the pushforward $(p_{\lambda})_{\ast}(\widetilde{\nu}_{\lambda,\mu})$ of the joint spectral measure $\widetilde{\nu}_{\lambda,\mu}$ under the orthogonal projection $p_{\lambda}$ onto the spectrum $\sigma(\lambda)$. See \cite[Section 2.5]{evans/pugh:2012iv} for more details.

The following $S_3$-invariant measures on $\mathbb{T}^2$, defined in \cite[Definition 1]{evans/pugh:2010i}, will be useful later. Note that these measures are also invariant under $D_{12}$.
\begin{Def} \label{def:4measures}
Let $\omega = e^{2 \pi i/3}$, $\tau = e^{2 \pi i/n}$. We define the following measures on $\mathbb{T}^2$:
\begin{enumerate}
\item $\mathrm{d}^{((n))}$, the uniform Dirac measure on the $S_3$-orbit of the points $(\tau, \tau)$, $(\overline{\omega} \, \overline{\tau}, \omega)$, $(\omega, \overline{\omega} \, \overline{\tau})$, for $n \in \mathbb{Q}$, $n \geq 2$.
\item $\mathrm{d}^{(n,k)}$, the uniform Dirac measure on the $S_3$-orbit of the points $(\tau \, e^{2 \pi i k}, \tau)$, $(\tau, \tau \, e^{2 \pi i k})$, $(\overline{\omega} \, \overline{\tau}, \omega \, e^{2 \pi i k})$, $(\omega \, e^{2 \pi i k}, \overline{\omega} \, \overline{\tau})$, $(\overline{\omega} \, \overline{\tau} \, e^{-2 \pi i k}, \omega \, e^{-2 \pi i k})$, $(\omega \, e^{-2 \pi i k}, \overline{\omega} \, \overline{\tau} \, e^{-2 \pi i k})$, for $n,k \in \mathbb{Q}$, $n > 2$, $0 \leq k \leq 1/n$.
\end{enumerate}
\end{Def}

The sets $\mathrm{Supp}(\mathrm{d}^{((n))})$, $\mathrm{Supp}(\mathrm{d}^{(n,k)})$ are illustrated in \cite[Figures 4, 5]{evans/pugh:2010i}.
For $n > 2$ and $0 < k < 1/n$, $|\mathrm{Supp}(\mathrm{d}^{((n))})| = 18$, whilst $|\mathrm{Supp}(\mathrm{d}^{(n,k)})| = 36$. The cardinalities of the other sets are $|\mathrm{Supp}(\mathrm{d}^{(n,0)})| = |\mathrm{Supp}(\mathrm{d}^{(n,1/n)})| = 18$ for $n > 2$, and $|\mathrm{Supp}(\mathrm{d}^{((2))})| = 9$.
Some relations between these measures are given in \cite[Section 2]{evans/pugh:2010i}.

\section{Spectral measures for ${}^W \hspace{-2mm} \mathcal{A}_{\infty}(G_2)$} \label{sect:measure:A12infty(G2)}

For the remainder of the paper we focus on the fundamental representations $\rho_1$ and $\rho_2$ of $G_2$. We first consider their restrictions to $\mathbb{T}^2$.
As discussed in Section \ref{sect:rep_theoryG2T2}, their corresponding McKay graphs are ${}^W \hspace{-2mm} \mathcal{A}^{\rho_j}_{\infty}(G_2)$, illustrated in Figures \ref{fig-A_12infty(G2)1}, \ref{fig-A_12infty(G2)2} for $j=1,2$ respectively.

We define commuting self-adjoint operators which may be identified with the adjacency matrix of ${}^W \hspace{-2mm} \mathcal{A}^{\rho_j}_{\infty}(G_2)$. We define operators $v_Z^1$, $v_Z^2$ on $\ell^2(\mathbb{Z}) \otimes \ell^2(\mathbb{Z})$ by
\begin{align}
v_Z^1 & = 1 \otimes 1 + s \otimes 1 + s^{\ast} \otimes 1 + 1 \otimes s + 1 \otimes s^{\ast} + s \otimes s^{\ast} + s^{\ast} \otimes s, \\
v_Z^2 & = v_Z^1 + 1 \otimes 1 + s \otimes s + s^{\ast} \otimes s^{\ast} + s^2 \otimes s^{\ast} + (s^{\ast})^2 \otimes s + s \otimes (s^{\ast})^2 + s^{\ast} \otimes s^2, \qquad
\end{align}
where the unitary $s$ is the bilateral shift on $\ell^2(\mathbb{Z})$.
Let $\Omega$ denote the vector $(\delta_{i,0})_i$. Then $v_Z^j$ is identified with the adjacency matrix of ${}^W \hspace{-2mm} \mathcal{A}^{\rho_j}_{\infty}(G_2)$, $j=1,2$, where we regard the vector $\Omega \otimes \Omega$ as corresponding to the vertex $(0,0)$ of ${}^W \hspace{-2mm} \mathcal{A}^{\rho_j}_{\infty}(G_2)$, and the operators of the form $s^l \otimes s^m$ which appear as terms in $v_Z^j$ as corresponding to the edges on ${}^W \hspace{-2mm} \mathcal{A}^{\rho_j}_{\infty}(G_2)$. Then $(s^{\lambda_1} \otimes s^{\lambda_2})(\Omega \otimes \Omega)$ corresponds to the vertex $(\lambda_1,\lambda_2)$ of ${}^W \hspace{-2mm} \mathcal{A}^{\rho_j}_{\infty}(G_2)$ for any $\lambda_1,\lambda_2 \in \mathbb{Z}$, and applying $(v_Z^j)^m$ to $\Omega \otimes \Omega$ gives a vector $y=(y_{(\lambda_1,\lambda_2)})$ in $\ell^2({}^W \hspace{-2mm} \mathcal{A}^{\rho_j}_{\infty}(G_2))$, where $y_{(\lambda_1,\lambda_2)}$ gives the number of paths of length $m$ on ${}^W \hspace{-2mm} \mathcal{A}^{\rho_j}_{\infty}(G_2)$ from $(0,0)$ to the vertex $(\lambda_1,\lambda_2)$.

We define a state $\varphi$ on $C^{\ast}(v_Z^1,v_Z^2)$ by $\varphi( \, \cdot \, ) = \langle \, \cdot \, (\Omega \otimes \Omega), \Omega \otimes \Omega \rangle$.
We use the notation $(a_1,a_2,\ldots,a_k)!$ to denote the multinomial coefficient $(\sum_{i=1}^k a_i)!/\prod_{i=1}^k(a_i!)$.
Then we have cross moments
\begin{align}
\varsigma_{m,n} & = \varphi((v_Z^1)^m(v_Z^2)^n) \nonumber \\
& = \sum_{\stackrel{k_i \geq 0}{\scriptscriptstyle{\sum_i k_i \leq m}}} \sum_{\stackrel{l_i \geq 0}{\scriptscriptstyle{\sum_i l_i \leq m}}} (k_1, k_2, \ldots, k_6, m - \sum_i k_i)! (l_1, l_2, \ldots, l_{12}, n - \sum_i l_i)! \; \varphi(s^{r_1} \otimes s^{r_2}) \nonumber \\
& = \sum_{\stackrel{k_i \geq 0}{\scriptscriptstyle{\sum_i k_i \leq m}}} \sum_{\stackrel{l_i \geq 0}{\scriptscriptstyle{\sum_i l_i \leq m}}} (k_1, k_2, \ldots, k_6, m - \sum_i k_i)! (l_1, l_2, \ldots, l_{12}, n - \sum_i l_i)! \; \delta_{r_1, 0} \; \delta_{r_2, 0}, \label{eqn:crossmomentsA_12infty(G2)}
\end{align}
where
\begin{align}
r_1 & = k_1-k_2+k_5-k_6 + l_1-l_2+l_5-l_6+l_7-l_8+2l_9-2l_{10}+l_{11}-l_{12}, \label{eqn:r1-G2} \\
r_2 & = k_3-k_4-k_5+k_6 + l_3-l_4-l_5+l_6+l_7-l_8-l_9+l_{10}-2l_{11}+2l_{12}, \label{eqn:r2-G2}
\end{align}
If $n=0$ then $l_i=0$ for all $i$, and we get a non-zero contribution when $k_4 = k_1-k_2+k_3$ and $k_6 = k_1-k_2+k_5$. So we obtain
\begin{equation} \label{eqn:momentsA_12infty(G2)1}
\varphi((v_Z^1)^m) = \sum_{k_i} (k_1, k_2, k_3, k_1-k_2+k_3, k_5, k_1-k_2+k_5, m-3k_1+k_2-2k_3-2k_5)!
\end{equation}
where the summation is over all integers $k_1,k_2,k_3,k_5 \geq 0$ such that $0 \leq k_1-k_2+k_3, k_1-k_2+k_5, 3k_1-k_2+2k_3+2k_5 \leq m$.
If $m=0$ then $k_i=0$ for all $i$, and we get a non-zero contribution when $l_4 = l_1-l_2+l_3+2l_7-2l_8+l_9-l_{10}-l_{11}+l_{12}$ and $l_6 = l_1-l_2+l_5+l_7-l_8+2l_9-2l_{10}+l_{11}-l_{12}$. So we obtain
\begin{equation} \label{eqn:momentsA_12infty(G2)2}
\varphi((v_Z^2)^n) = \sum_{l_i} 2^{n-p_3} (l_1, l_2, l_3, p_1, l_5, p_2, l_7, l_8, l_9, l_{10}, l_{11}, l_{12}, n-p_3)!
\end{equation}
where $p_1 = l_1-l_2+l_3+2l_7-2l_8+l_9-l_{10}-l_{11}+l_{12}$, $p_2 = l_1-l_2+l_5+l_7-l_8+2l_9-2l_{10}+l_{11}-l_{12}$, $p_3 = 3l_1-l_2+2l_3+2l_5+4l_7-2l_8+4l_9-2l_{10}+l_{11}+l_{12}$ and the summation is over all integers $l_1,l_2,l_3,l_5,l_7,l_8,l_9,l_{10},l_{11},l_{12} \geq 0$ such that $0 \leq p_1, p_2, p_3 \leq n$.

\subsection{Joint spectral measure for ${}^W \hspace{-2mm} \mathcal{A}_{\infty}(G_2)$ over $\mathbb{T}^2$}

The ranges of the restrictions (\ref{eqn:restriction-chi1}), (\ref{eqn:restriction-chi2}) of the characters $\chi_j$, $j=1,2$, of the fundamental representations of $G_2$ to $\mathbb{T}^2$ are given by $I_j$, where $I_1 := I_{\rho_1} = \{ 1 + 2\mathrm{Re}(\omega_1) + 2\mathrm{Re}(\omega_2) + 2\mathrm{Re}(\omega_1\omega_2^{-1}) | \, \omega_1,\omega_2 \in \mathbb{T} \} = [-2,7]$ and $I_2 := I_{\rho_2} = \{ 2 + 2\mathrm{Re}(\omega_1) + 2\mathrm{Re}(\omega_2) + 2\mathrm{Re}(\omega_1\omega_2^{-1}) + 2\mathrm{Re}(\omega_1\omega_2) + 2\mathrm{Re}(\omega_1^2\omega_2^{-1}) + 2\mathrm{Re}(\omega_1\omega_2^{-2}) | \, \omega_1,\omega_2 \in \mathbb{T} \} = [-2,14]$.
Since the spectrum $\sigma(s)$ of $s$ is $\mathbb{T}$, the spectrum $\sigma(v_Z^1)$ of $v_Z^1$ is $I_1 = [-2,7]$, and the spectrum $\sigma(v_Z^2)$ of $v_Z^2$ is $I_2 = [-2,14]$.
We now determine the $D_{12}$-invariant spectral measure $\varepsilon$ on $\mathbb{T}^2$ for the graphs ${}^W \hspace{-2mm} \mathcal{A}^{\rho_j}_{\infty}(G_2)$, $j=1,2$.

\begin{Thm} \label{thm:measureG2-T2}
The joint spectral measure $\varepsilon(\omega_1,\omega_2)$ (on $\mathbb{T}^2$) for the graphs ${}^W \hspace{-2mm} \mathcal{A}^{\rho_j}_{\infty}(G_2)$, $j=1,2$, is given by the uniform Lebesgue measure $\mathrm{d}\varepsilon(\omega_1,\omega_2) = \mathrm{d}\omega_1 \, \mathrm{d}\omega_2$.
\end{Thm}
\emph{Proof:}
The $m,n^{\mathrm{th}}$ cross moment is given by
\begin{eqnarray*}
\lefteqn{ \int_{\mathbb{T}^2} (\chi_1(\omega_1, \omega_2))^m (\chi_2(\omega_1, \omega_2))^n \mathrm{d}\omega_1 \, \mathrm{d}\omega_2 } \\
& = & \sum_{\stackrel{k_i \geq 0}{\scriptscriptstyle{\sum_i k_i \leq m}}} \sum_{\stackrel{l_i \geq 0}{\scriptscriptstyle{\sum_i l_i \leq m}}} (k_1, k_2, \ldots, k_6, m - \sum_i k_i)! (l_1, l_2, \ldots, l_{12}, n - \sum_i l_i)! \; \int_{\mathbb{T}^2} \omega_1^{r_1} \omega_2^{r_2} \mathrm{d}\omega_1 \, \mathrm{d}\omega_2 \\
& = & \sum_{\stackrel{k_i \geq 0}{\scriptscriptstyle{\sum_i k_i \leq m}}} \sum_{\stackrel{l_i \geq 0}{\scriptscriptstyle{\sum_i l_i \leq m}}} (k_1, k_2, \ldots, k_6, m - \sum_i k_i)! (l_1, l_2, \ldots, l_{12}, n - \sum_i l_i)! \; \delta_{r_1, 0} \; \delta_{r_2, 0},
\end{eqnarray*}
where $r_1$, $r_2$ are as in (\ref{eqn:r1-G2}), (\ref{eqn:r2-G2}), since $\int_{\mathbb{T}} u^m \mathrm{d}u = \delta_{m,0}$. This is equal to the cross moments $\varphi((v_Z^1)^m(v_Z^2)^n)$ given in (\ref{eqn:crossmomentsA_12infty(G2)}).
\hfill
$\Box$

In fact, the measure $\varepsilon(\omega_1,\omega_2)$ given above is the joint spectral measure over $\mathbb{T}^2$ for the pair of McKay graphs (${}^W \hspace{-2mm} \mathcal{A}_{\infty}^{\lambda}(G_2)$, ${}^W \hspace{-2mm} \mathcal{A}_{\infty}^{\mu}(G_2)$) for any pair $\lambda,\mu$ of irreducible representations of $G_2$, by a similar proof. Thus the spectral measure over $\mathbb{T}^2$ is independent of the choice of irreducible representations used to construct the McKay graphs.

\subsection{Joint spectral measure for ${}^W \hspace{-2mm} \mathcal{A}_{\infty}(G_2)$ on $\mathfrak{D}$} \label{sect:measureA12inftyG2-D}

Let $x := x_{(1,0)} = \chi_1(\omega_1,\omega_2)$ and $y := x_{(1,1)} = \chi_2(\omega_1,\omega_2)$, or explicitly,
\begin{align}
x & = 1 + 2\cos(2\pi\theta_1) + 2\cos(2\pi\theta_2) + 2\cos(2\pi(\theta_1-\theta_2)), \label{eqn:x-G2} \\
y & = x + 1 + 2\cos(2\pi(\theta_1+\theta_2)) + 2\cos(2\pi(2\theta_1-\theta_2)) + 2\cos(2\pi(\theta_1-2\theta_2)), \qquad \label{eqn:y-G2}
\end{align}
and denote by $\Psi$ the map $\Psi_{(1,0),(1,1)}: (\omega_1,\omega_2) \mapsto (x,y)$.

\begin{figure}[tb]
\begin{center}
  \includegraphics[width=60mm]{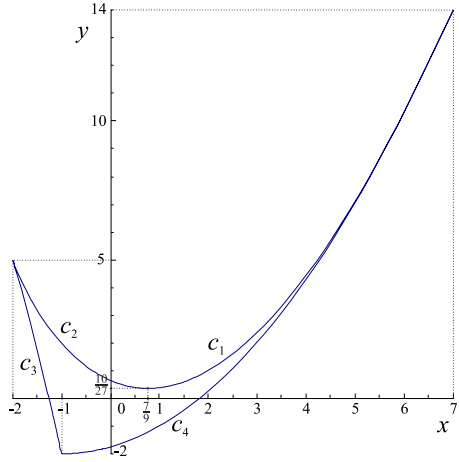}\\
 \caption{The domain $\mathfrak{D} = \Psi(F)$ for $G_2$.} \label{fig:DomainD-G2}
\end{center}
\end{figure}

We now describe $\mathfrak{D} := \mathfrak{D}_{(1,0),(1,1)}$, illustrated in Figure \ref{fig:DomainD-G2}, which is the joint spectrum $\sigma(v_Z^1,v_Z^2)$ of the commuting self-adjoint operators $v_Z^1$, $v_Z^2$.
The boundary of $F$ given by $\theta_1 = 2\theta_2$ yields the curves $c_1$, $c_2$ for $\theta_2 \in [0,1/3]$, which are both given by the parametric equations
\begin{align*}
x & = 1 + 4\cos(2\pi\theta_2) + 2\cos(4\pi\theta_2) = -1 + 4\cos(2\pi\theta_2) + 4\cos^2(2\pi\theta_2), \\
y & = 4 + 4\cos(2\pi\theta_2) + 2\cos(4\pi\theta_2) + 4\cos(6\pi\theta_2) \\
& = 2 - 8\cos(2\pi\theta_2) + 4\cos^2(2\pi\theta_2) + 16\cos^3(2\pi\theta_2).
\end{align*}
The boundary of $F$ given by $\theta_1 = -\theta_2$ yields the curve $c_3$ given by the parametric equations
\begin{align*}
x & = -1 + 4\cos(2\pi\theta_1) + 4\cos^2(2\pi\theta_1), \\
y & = 2 - 8\cos(2\pi\theta_1) + 4\cos^2(2\pi\theta_1) + 16\cos^3(2\pi\theta_1),
\end{align*}
where $\theta_1 \in [1/2,2/3]$.
Finally, the boundary $\theta_1 = \theta_2$ of $F$ yields the curve $c_4$ given by the parametric equations
\begin{align*}
x & = 3 + 4\cos(2\pi\theta_1), \\
y & = 4 + 8\cos(2\pi\theta_1) + 2\cos(4\pi\theta_1) = 2 + 8\cos(2\pi\theta_1) + 4\cos^2(2\pi\theta_1),
\end{align*}
where $\theta_1 \in [1/2,1]$.
As functions of $x$, the boundaries of $\mathfrak{D}$ are obtained by writing $\cos(2\pi \theta)$ in terms of $x$ in the above parametric equations, which are at worst quadratic in $\cos(2\pi \theta)$. The boundaries of $\mathfrak{D}$ are thus given by the curves \cite{uhlmann/meinel/wipf:2007}
\begin{align}
&c_1: & \hspace{-10mm} y & = -5(x+1)+2(x+2)^{3/2}, \qquad x \in [-2,7/9], \label{eqn:boundaryD-y-G2-1} \\
&c_2: & \hspace{-10mm} y & = -5(x+1)+2(x+2)^{3/2}, \qquad x \in [7/9,7], \label{eqn:boundaryD-y-G2-1a} \\
&c_3: & \hspace{-10mm} y & = -5(x+1)-2(x+2)^{3/2}, \qquad x \in [-2,-1], \label{eqn:boundaryD-y-G2-2}\\
&c_4: & \hspace{-10mm} 4y & = x^2+2x-7, \hspace{30mm} x \in [-1,7]. \label{eqn:boundaryD-y-G2-3}
\end{align}
On the other hand, writing these curves as function of $y$ involves cubic equations in $\cos(2\pi \theta)$, and the boundaries of $\mathfrak{D}$ are given by the curves
\begin{align}
&c_1: & \hspace{-10mm} x & = -1+4p_2(y)+4p_2(y)^2, \qquad y \in [10/27,14], \label{eqn:boundaryD-x-G2-1} \\
&c_2: & \hspace{-10mm} x & = -1+4p_3(y)+4p_3(y)^2, \qquad y \in [10/27,5] \label{eqn:boundaryD-x-G2-1a} \\
&c_3: & \hspace{-10mm} x & = -1+4p_1(y)+4p_1(y)^2, \qquad y \in [-2,5], \label{eqn:boundaryD-x-G2-2} \\
&c_4: & \hspace{-10mm} x & = -1+2(y+2)^{1/2}, \hspace{18mm} y \in [-2,14], \label{eqn:boundaryD-x-G2-3}
\end{align}
where $p_i$ is given by $12p_i(y) = -1 - \epsilon_i P - 25\overline{\epsilon_i} P^{-1}$, for $\epsilon_j = e^{2 \pi i (j-1)/3}$ and $P = (145-54y + 2\sqrt{3^3(27y^2-145y+50)})^{1/3}$, and we take the positive square root in equation (\ref{eqn:boundaryD-x-G2-3}).

\begin{figure}[tb]
\begin{minipage}[t]{8.9cm}
\begin{center}
  \includegraphics[width=80mm]{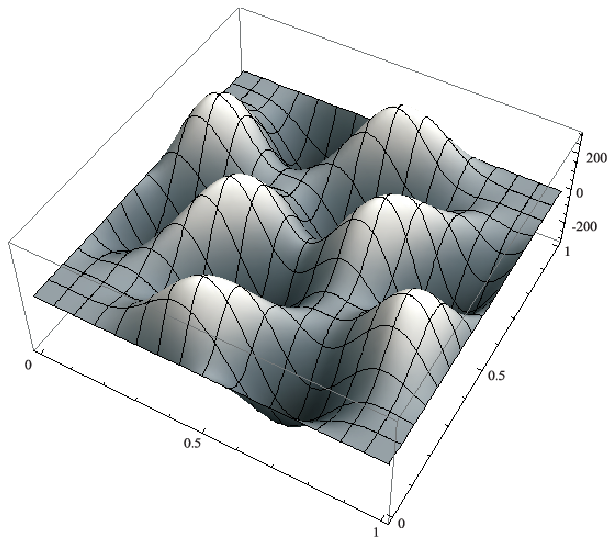}\\
 \caption{The Jacobian $J$ over $\mathbb{T}^2$.} \label{fig:J_overT2-3Dplot-G2}
\end{center}
\end{minipage}
\hfill
\begin{minipage}[t]{6.9cm}
\begin{center}
  \includegraphics[width=55mm]{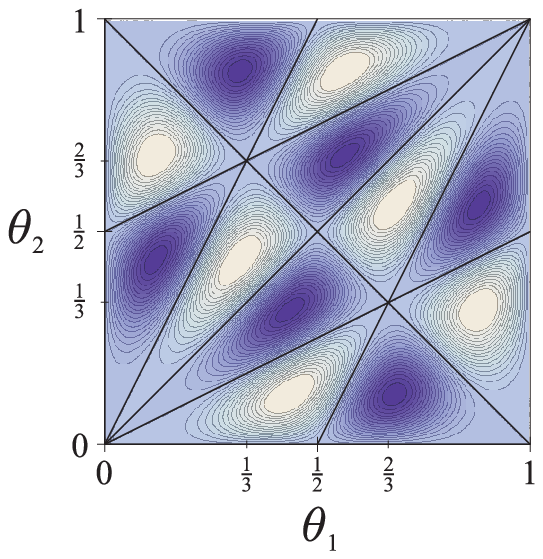}\\
 \caption{Contour plot of $J$ over $\mathbb{T}^2$.} \label{fig:J_overT2-contours-G2}
\end{center}
\end{minipage}
\end{figure}

Under the change of variables (\ref{eqn:x-G2}), (\ref{eqn:y-G2}), the Jacobian $J = \mathrm{det}(\partial(x,y)/\partial(\theta_1,\theta_2))$ is given by
\begin{align} \label{eqn:J[theta]-G2}
J (\theta_1,\theta_2) & = 8 \pi^2 (\cos(2 \pi (2\theta_1 + \theta_2)) + \cos(2 \pi (\theta_1 - 3\theta_2)) + \cos(2 \pi (3\theta_1 - 2\theta_2)) \nonumber \\
& \qquad - \cos(2 \pi (\theta_1 + 2\theta_2)) - \cos(2 \pi (3\theta_1 - \theta_2)) - \cos(2 \pi (2\theta_1 - 3\theta_2))). \qquad
\end{align}
The Jacobian is real and is illustrated in Figures \ref{fig:J_overT2-3Dplot-G2}, \ref{fig:J_overT2-contours-G2}, where its values are plotted over the torus $\mathbb{T}^2$.

With $\omega_j = e^{2 \pi i \theta_j}$, $j=1,2$, the Jacobian is given in terms of $\omega_1, \omega_2 \in \mathbb{T}$ by
\begin{align}
J (\omega_1,\omega_2) & = 8 \pi^2 \mathrm{Re}(\omega_1^2\omega_2 + \omega_1\omega_2^{-3} + \omega_1^3\omega_2^{-2} - \omega_1\omega_2^2 - \omega_1^3\omega_2^{-1} - \omega_1^2\omega_2^{-3}) \nonumber \\
& = 4 \pi^2 (\omega_1^2\omega_2 + \omega_1^{-2}\omega_2^{-1} + \omega_1\omega_2^{-3} + \omega_1^{-1}\omega_2^{3} + \omega_1^3\omega_2^{-2} + \omega_1^{-3}\omega_2^{2} \nonumber \\
& \qquad - \omega_1\omega_2^2 - \omega_1^{-1}\omega_2^{-2} - \omega_1^3\omega_2^{-1} - \omega_1^{-3}\omega_2 - \omega_1^2\omega_2^{-3} - \omega_1^{-2}\omega_2^{3}). \qquad \label{eqn:J[omega]-G2}
\end{align}
The Jacobian $J$ is invariant under $T_6 \in D_{12}$, whilst $T_2(J) = -J$. Thus $J^2$ is invariant under the action of $D_{12}$, and we seek an expression for $J^2$ in terms of the $D_{12}$-invariant variables $x$, $y$, which may be obtained as a product of the roots appearing as the equations of the boundary of $\mathfrak{D}$ in (\ref{eqn:boundaryD-y-G2-1})-(\ref{eqn:boundaryD-y-G2-3}), and is given (up to a factor of $16 \pi^4$) as (see also \cite{uhlmann/meinel/wipf:2007})
\begin{equation} \label{J2-G2}
J^2(x,y) = (4x^3-x^2-2x-10xy-y^2-10y+7)(x^2+2x-7-4y),
\end{equation}
for $(x,y) \in \mathfrak{D}$, which can easily be checked by substituting for $x$, $y$ as in (\ref{eqn:x-G2}), (\ref{eqn:y-G2}).
Note that the Jacobian (\ref{J2-G2}) is a cubic in $y$, with the three roots appearing as the equations of the boundary of $\mathfrak{D}$ in (\ref{eqn:boundaryD-y-G2-1})-(\ref{eqn:boundaryD-y-G2-3}).
However, although the Jacobian (\ref{J2-G2}) is a quintic in $x$, only four of the roots appear as the equations of the boundary of $\mathfrak{D}$ in (\ref{eqn:boundaryD-x-G2-1})-(\ref{eqn:boundaryD-x-G2-3}). The fifth root $x = -1-2(y+2)^{1/2}$ only intersects with $\mathfrak{D}$ at the point $(-1,-2)$.
The factorization of $J$ in (\ref{J2-G2}) and the equations for the boundaries of $\mathfrak{D}$ given in (\ref{eqn:boundaryD-y-G2-1})-(\ref{eqn:boundaryD-x-G2-3}) will be used in Sections \ref{sect:measureA12inftyG2-R}, \ref{sect:measureAinftyG2-R} to determine explicit expressions for the weights which appear in the spectral measures $\mu_{v_Z^j}$ over $I_j$ in terms of elliptic integrals.
From (\ref{J2-G2}) we see that the Jacobian vanishes only on the boundary of $\mathfrak{D}$, or over $\mathbb{T}^2$ only on the boundaries of the images of the fundamental domain $F$ under $D_{12}$.

Since $J$ is real, $J^2 \geq 0$, and we have the following expressions for the Jacobian $J$:
\begin{align*}
J (\theta_1,\theta_2) & = 8 \pi^2 (\cos(2 \pi (2\theta_1 + \theta_2)) + \cos(2 \pi (\theta_1 - 3\theta_2)) + \cos(2 \pi (3\theta_1 - 2\theta_2)) \\
& \qquad - \cos(2 \pi (\theta_1 + 2\theta_2)) - \cos(2 \pi (3\theta_1 - \theta_2)) - \cos(2 \pi (2\theta_1 - 3\theta_2))), \\
J (\omega_1,\omega_2) & = 4 \pi^2 (\omega_1^2\omega_2 + \omega_1^{-2}\omega_2^{-1} + \omega_1\omega_2^{-3} + \omega_1^{-1}\omega_2^{3} + \omega_1^3\omega_2^{-2} + \omega_1^{-3}\omega_2^{2} \\
& \qquad - \omega_1\omega_2^2 - \omega_1^{-1}\omega_2^{-2} - \omega_1^3\omega_2^{-1} - \omega_1^{-3}\omega_2 - \omega_1^2\omega_2^{-3} - \omega_1^{-2}\omega_2^{3}), \\
|J (x,y)| & = 4 \pi^2 \sqrt{(4x^3-x^2-2x-10xy-y^2-10y+7)(x^2+2x-7-4y)},
\end{align*}
where $0 \leq \theta_1,\theta_2 < 1$, $\omega_1,\omega_2 \in \mathbb{T}$ and $(x,y) \in \mathfrak{D}$. Note that the expression under the square root is real and non-negative since $J^2$ is.

Then
\begin{equation} \label{eqn:integral_overDj-G2}
\int_{F} \psi(\chi_1(\omega_1,\omega_2),\chi_2(\omega_1,\omega_2)) \mathrm{d}\omega_1 \, \mathrm{d}\omega_2 = \int_{\mathfrak{D}} \psi(x,y) |J(x,y)|^{-1} \mathrm{d}x \, \mathrm{d}y,
\end{equation}
and from Theorem \ref{thm:measureG2-T2} and (\ref{eqn:measureT2=12C}) we obtain
\begin{Thm}
The joint spectral measure $\widetilde{\nu}$ (over $\mathfrak{D}$) for the graphs ${}^W \hspace{-2mm} \mathcal{A}^{\rho_j}_{\infty}(G_2)$, $j=1,2$, is
$$\mathrm{d}\widetilde{\nu}(x,y) = 12 \, |J(x,y)|^{-1} \mathrm{d}x \, \mathrm{d}y.$$
\end{Thm}

\subsection{Spectral measure for ${}^W \hspace{-2mm} \mathcal{A}_{\infty}(G_2)$ on $\mathbb{R}$} \label{sect:measureA12inftyG2-R}

We now compute the spectral measure $\nu_{v_Z^j} = \nu_{\rho_j}$ over $I_j$, which is determined by its moments $\varphi((v_Z^j)^m) = \int_{I_j} x_j^m \mathrm{d}\nu_{v_Z^j}(x_j)$ for all $m \in \mathbb{N}$, where $x_j = x,y$ for $j=1,2$ respectively.

\begin{figure}[tb]
\begin{minipage}[t]{5.5cm}
\begin{center}
  \includegraphics[width=51mm]{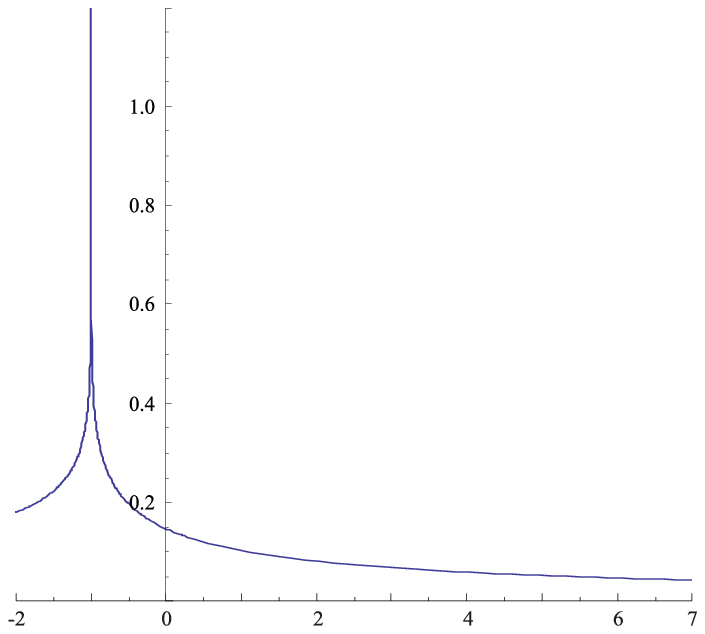}
 \caption{$J_1^{\mathbb{T}^2}(x)$} \label{fig-Jx-T2-G2}
\end{center}
\end{minipage}
\hfill
\begin{minipage}[t]{9.5cm}
\begin{center}
  \includegraphics[width=90mm]{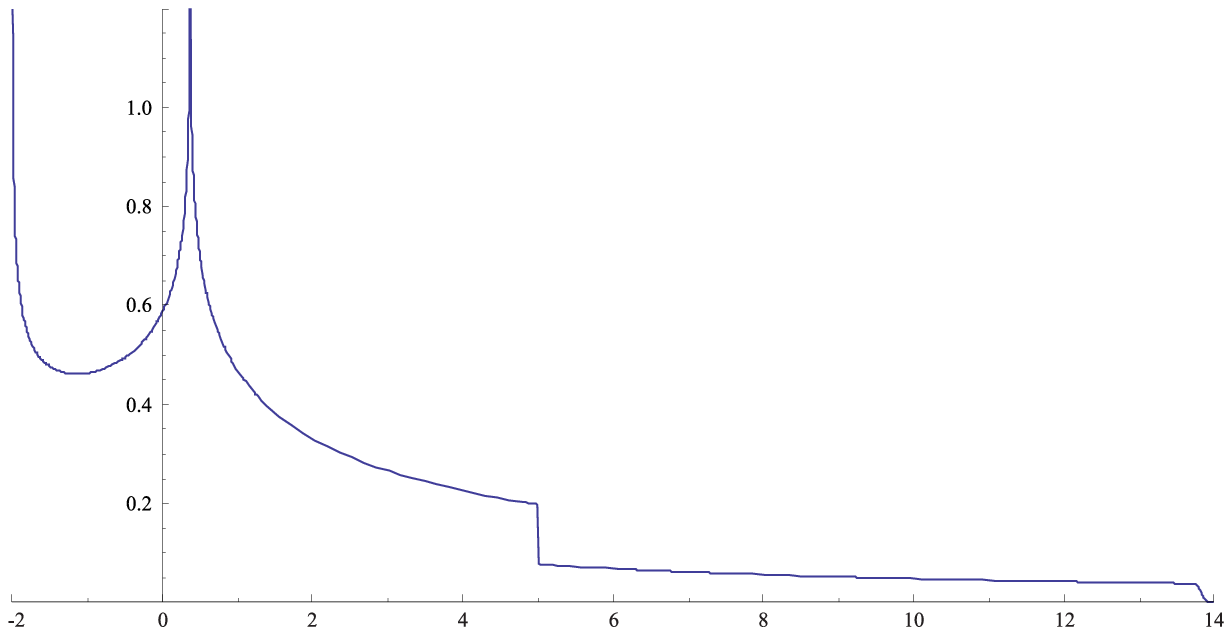}
 \caption{$J_2^{\mathbb{T}^2}(y)$} \label{fig-Jy-T2-G2}
\end{center}
\end{minipage}
\end{figure}

For $\nu_{v_Z^1}$ we set $\psi(x,y) = x^m$ in (\ref{eqn:integral_overDj-G2}) and integrate with respect to $y$. Similarly, setting $\psi(x,y) = y^m$ in (\ref{eqn:integral_overDj-G2}), the measure $\nu_{v_Z^2}$ is obtained by integrating with respect to $x$. More explicitly, using the expressions for the boundaries of $\mathfrak{D}$ given in (\ref{eqn:boundaryD-y-G2-1})-(\ref{eqn:boundaryD-x-G2-3}),
the spectral measure $\nu_{v_Z^1}$ (over $[-2,7]$) for the graph ${}^W \hspace{-2mm} \mathcal{A}^{\rho_1}_{\infty}(G_2)$ is $\mathrm{d}\nu_{v_Z^1}(x) = J_1^{\mathbb{T}^2}(x) \, \mathrm{d}x$, where $J_1^{\mathbb{T}^2}(x)$ is given by
$$J_1^{\mathbb{T}^2}(x) = \left\{
\begin{array}{cl}
\displaystyle{12 \int_{-5(x+1)-2(x+2)^{3/2}}^{-5(x+1)+2(x+2)^{3/2}} |J(x,y)|^{-1} \, \mathrm{d}y} & \textrm{ for } x \in [-2,-1], \\
\displaystyle{12 \int_{(x^2+2x-7)/4}^{-5(x+1)+2(x+2)^{3/2}} |J(x,y)|^{-1} \, \mathrm{d}y} & \textrm{ for } x \in [-1,7].
\end{array} \right. $$
The weight $J_1^{\mathbb{T}^2}(x)$ is an integral of the reciprocal of the square root of a cubic in $y$, and thus can be written in terms of the complete elliptic integral $K(m)$ of the first kind, $K(m) = \int_0^{\pi/2} (1-m\sin^2\theta)^{-1/2} \mathrm{d}\theta$. Using \cite[Eqn. 235.00]{byrd/friedman:1971},
$$\frac{6}{\pi^2 \left( 8(x+2)^{3/2} - x^2-22x-13 \right)^{1/2}} \; K(v(x)) = \frac{3 \, v(x)^{1/2}}{2\pi^2 (x+2)^{3/4}} \; K(v(x))$$
for $x \in [-2,-1]$, where $v(x) = 16(x+2)^{3/2}/(8(x+2)^{3/2}-x^2-22x-13)$, whilst for $x \in [-1,7]$,
$$\frac{6 \, v(x)^{-1/2}}{\pi^2 \left( 8(x+2)^{3/2} - x^2-22x-13 \right)^{1/2}} \; K(v(x)^{-1}) = \frac{3}{2\pi^2 (x+2)^{3/4}} \; K(v(x)^{-1}).$$
The weight $J_1^{\mathbb{T}^2}(x)$ is illustrated in Figure \ref{fig-Jx-T2-G2}, up to a factor $4\pi^2$.

The spectral measure $\nu_{v_Z^2}$ (over $[-2,14]$) for the graph ${}^W \hspace{-2mm} \mathcal{A}^{\rho_2}_{\infty}(G_2)$ is $\mathrm{d}\nu_{v_Z^2}(y) = J_2^{\mathbb{T}^2}(y) \, \mathrm{d}y$, where $J_2^{\mathbb{T}^2}(y)$ is given by
\begin{eqnarray*}
& \displaystyle{12 \int_{-1+4p_1(y)+4p_1(y)^2}^{-1+2(y+2)^{1/2}} |J(x,y)|^{-1} \, \mathrm{d}x} & \textrm{ for } y \in [-2,10/27], \\
& \displaystyle{12 \int_{-1+4p_1(y)+4p_1(y)^2}^{-1+4p_3(y)+4p_3(y)^2} |J(x,y)|^{-1} \, \mathrm{d}x + 12 \int_{-1+4p_2(y)+4p_2(y)^2}^{-1+2(y+2)^{1/2}} |J(x,y)|^{-1} \, \mathrm{d}x} & \textrm{ for } y \in [10/27,5], \\
& \displaystyle{12 \int_{-1+4p_2(y)+4p_2(y)^2}^{-1+2(y+2)^{1/2}} |J(x,y)|^{-1} \, \mathrm{d}x} & \textrm{ for } y \in [5,14].
\end{eqnarray*}
A numerical plot of the weight $J_2^{\mathbb{T}^2}(y)$ is illustrated in Figure \ref{fig-Jy-T2-G2}, again up to a factor $4\pi^2$.

\section{Spectral measures for $\mathcal{A}_{\infty}(G_2)$} \label{sect:measureAinftyG2}

We now consider the fundamental representations $\rho_j$, $j=1,2$, of $G_2$.
As discussed in Section \ref{sect:rep_theoryG2T2}, their corresponding McKay graphs are $\mathcal{A}^{\rho_j}_{\infty}(G_2)$, $j=1,2$.
This section follows the same arguments as \cite[$\S$6.2]{evans/pugh:2009v}. The new feature here is the presence of terms such as $l^{m}(l^{\ast})^{p}$ for $m \leq p$, where $l$ is the unilateral shift to the right on $\ell^2(\mathbb{N})$, which correspond to the fact that certain edges on the graphs $\mathcal{A}_{\infty}^{\rho_j}(G_2)$, $j=1,2$, only appear when far enough away from the boundary of the graph.

We define self-adjoint operators $v_N^1$, $v_N^2$ on $\ell^2(\mathbb{N}) \otimes \ell^2(\mathbb{N})$ by
\begin{align}
v_N^1 & = ll^{\ast} \otimes 1 + l \otimes 1 + l^{\ast} \otimes 1 + l^{\ast} \otimes l + l \otimes l^{\ast} + l^2 \otimes l^{\ast} + (l^{\ast})^2 \otimes l, \label{eqn:vN1} \\
v_N^2 & = ll^{\ast} \otimes 1 + l^2 l^{\ast} \otimes 1 + l(l^{\ast})^2 \otimes 1 + l(l^{\ast})^2 \otimes l + l^2 l^{\ast} \otimes l^{\ast} + l^2 \otimes l^{\ast} + (l^{\ast})^2 \otimes l \qquad \nonumber \\
& \quad + 1 \otimes ll^{\ast} + 1 \otimes l + 1 \otimes l^{\ast} + l^3 \otimes l^{\ast} + (l^{\ast})^3 \otimes l + l^3 \otimes (l^{\ast})^2 + (l^{\ast})^3 \otimes l^2. \label{eqn:vN2}
\end{align}
Let $\Omega$ denote the vector $(\delta_{i,0})_i$. Then $v_N^j$ is identified with the adjacency matrix of $\mathcal{A}^{\rho_j}_{\infty}(G_2)$, $j=1,2$, where we regard the vector $\Omega \otimes \Omega$ as corresponding to the vertex $(0,0)$ of $\mathcal{A}^{\rho_j}_{\infty}(G_2)$, and the operators of the form $l^{m_1}(l^{\ast})^{p_1} \otimes l^{m_2}(l^{\ast})^{p_2}$ which appear as terms in $v_N^j$ as corresponding to the edges on $\mathcal{A}^{\rho_j}_{\infty}(G_2)$. Then $(l^{\lambda_1} \otimes l^{\lambda_2})(\Omega \otimes \Omega)$ corresponds to the vertex $(\lambda_1,\lambda_2)$ of $\mathcal{A}^{\rho_j}_{\infty}(G_2)$ for any $\lambda_1,\lambda_2 \in \mathbb{N}$, and applying $(v_N^j)^m$ to $\Omega \otimes \Omega$ gives a vector $y=(y_{(\lambda_1,\lambda_2)})$ in $\ell^2(\mathcal{A}^{\rho_j}_{\infty}(G_2))$, where $y_{(\lambda_1,\lambda_2)}$ gives the number of paths of length $m$ on $\mathcal{A}^{\rho_j}_{\infty}(G_2)$ from $(0,0)$ to the vertex $(\lambda_1,\lambda_2)$. A term of the form $(l^{m}(l^{\ast})^{p} \otimes 1)(\lambda_1,\lambda_2)$, $m<p$, will be zero if $\lambda_1 <p$. Thus for example, $(ll^{\ast} \otimes 1)(0,\lambda_2) = 0$, which corresponds to the fact that there is no self-loop at the point $(0,\lambda_2)$ for any $\lambda_2 \in \mathbb{N}$, whereas $(ll^{\ast} \otimes 1)(\lambda_1,\lambda_2) = (\lambda_1,\lambda_2)$ for all $\lambda_1 \neq 0$, corresponding to the fact that there is a self-loop at these points.

It is not immediately obvious that these two operators commute. However this can be easily seen from the fact that $\mathcal{A}^{\rho_j}_{\infty}(G_2)$ are the multiplication graphs for the characters $\chi_j$ of the fundamental representations $\rho_j$, $j=1,2$, of $G_2$, where the vertices of $\mathcal{A}^{\rho_j}_{\infty}(G_2)$ are labeled by the characters of the irreducible representations of $G_2$.

Any vector $l^{p_1} \Omega \otimes l^{p_2} \Omega \in \ell^2(\mathbb{N}) \otimes \ell^2(\mathbb{N})$ can be written as a linear combination of elements of the form $(v_N^1)^{m_1} (v_N^2)^{m_2} (\Omega \otimes \Omega)$. This is not obvious from the definition of the operators $v_N^j$ given in (\ref{eqn:vN1}), (\ref{eqn:vN2}). However this also can be seen from the fact that $\mathcal{A}^{\rho_j}_{\infty}(G_2)$ are the multiplication graphs for the characters of the fundamental representations of $G_2$, where $\Omega \otimes \Omega$ corresponds to the character $\chi_{(0,0)}$ of the trivial representation. The characters of all other irreducible representations can be written as a linear combination of products of the form $\chi_1^m \chi_2^n \chi_{(0,0)}$, which follows from \cite[Proposition 1]{nesterenko/patera/tereszkiewicz:2011} where $X_1 = \chi_1-\chi_{(0,0)}$ and $X_2 = \chi_2-\chi_1-\chi_{(0,0)}$.

Thus the vector $\Omega \otimes \Omega$ is cyclic in $\ell^2(\mathbb{N}) \otimes \ell^2(\mathbb{N})$, and we have $\overline{C^{\ast}(v_N^1, v_N^2) (\Omega \otimes \Omega)} = \ell^2(\mathbb{N}) \otimes \ell^2(\mathbb{N})$.
We define a state $\varphi$ on $C^{\ast}(v_N^1, v_N^2)$ by $\varphi( \, \cdot \, ) = \langle \, \cdot \, (\Omega \otimes \Omega), \Omega \otimes \Omega \rangle$.
Since $C^{\ast}(v_N^1, v_N^2)$ is abelian and $\Omega \otimes \Omega$ is cyclic, we have that $\varphi$ is a faithful state on $C^{\ast}(v_N^1, v_N^2)$.
Then by \cite[Remark 2.3.2]{dykema/voiculescu/nica:1992} the support of $\widetilde{\nu}_{v_N^1, v_N^2}$ is equal to the joint spectrum $\sigma(v_N^1, v_N^2)$ of $v_N^1$, $v_N^2$.

Recall the decompositions of $\chi_1 \chi_{\mu}$ and $\chi_2 \chi_{\mu}$ for the characters of $G_2$ given in Section \ref{sect:rep_theoryG2T2}. These can be written as 
$\chi_j \chi_{\mu} = \sum_{\nu} \Delta_{\rho_j}(\mu,\nu) \chi_{\nu}$,
where $\Delta_{\rho_j}$ is the adjacency matrix of $\mathcal{A}_{\infty}^{\rho_j}(G_2)$, the McKay graph for the fundamental representation $\rho_j$ of $G_2$.
This equation can be interpreted as meaning that $v_N^j$ (identified with the adjacency matrix $\Delta_{\rho_j}$ of $\mathcal{A}_{\infty}^{\rho_j}(G_2)$) has eigenvector $(\chi_{\nu}(\theta))_{\nu}$ for eigenvalue $\chi_j(\theta)$, $\theta \in [0,2\pi]^2$.
Thus the spectrum of $v_N^j$ is given by $\chi_j(\mathbb{T}^2)$, and the joint spectrum $\sigma(v_N^1, v_N^2)$ is $\mathfrak{D}$.
The moments $\varphi((v_N^j)^m)$ count the number of closed paths of length $m$ on the graph $\mathcal{A}^{\rho_j}_{\infty}(G_2)$ which start and end at the apex vertex $(0,0)$.

\subsection{Joint spectral measure for $\mathcal{A}_{\infty}(G_2)$ over $\mathbb{T}^2$} \label{sect:measureAinftyG2-T2}

We prove in Section \ref{sect:measures_AkG2} that the measure given by
$\mathrm{d}\varepsilon(\omega_1,\omega_2) = J(\omega_1,\omega_2)^2 \mathrm{d}\omega_1 \, \mathrm{d}\omega_2 / 192 \pi^4$
is the joint spectral measure over $\mathbb{T}^2$ of $v_N^j$, $j=1,2$, where $\mathrm{d}\omega_l$ is the uniform Lebesgue measure on $\mathbb{T}$, $l=1,2$.
In fact, the measure $\varepsilon(\omega_1,\omega_2)$ is the joint spectral measure over $\mathbb{T}^2$ for the pair of McKay graphs ($\mathcal{A}_{\infty}^{\lambda}(G_2)$,$\mathcal{A}_{\infty}^{\mu}(G_2)$) for any pair $\lambda,\mu$ of irreducible representations of $G_2$.
We see that $\varepsilon$ is (up to some scalar) the reduced Haar measure of $G_2$ (c.f. \cite[$\S$6.3]{uhlmann/meinel/wipf:2007}).

\subsection{Spectral measure for $\mathcal{A}_{\infty}(G_2)$ on $\mathbb{R}$} \label{sect:measureAinftyG2-R}

\begin{figure}[tb]
\begin{minipage}[t]{5.5cm}
\begin{center}
  \includegraphics[width=51mm]{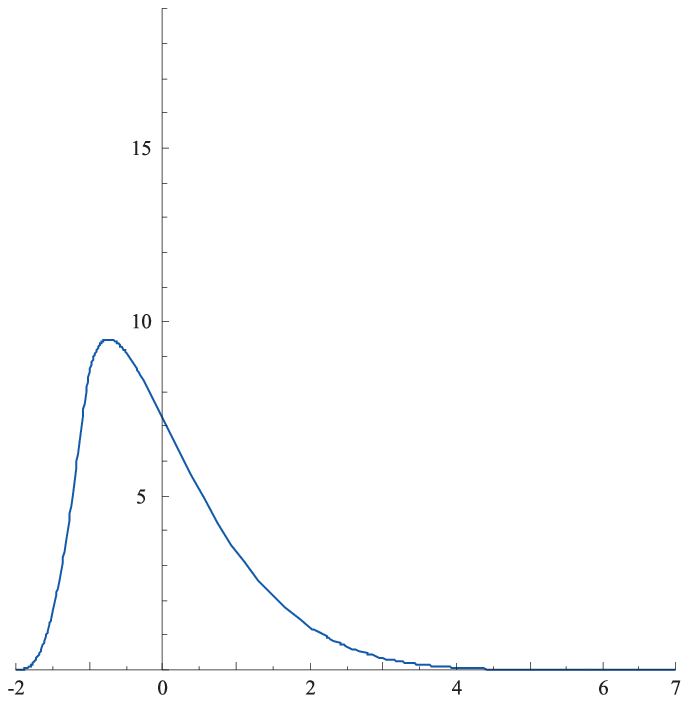}
 \caption{$J_1^{G_2}(x)$} \label{fig-Jx-G2}
\end{center}
\end{minipage}
\hfill
\begin{minipage}[t]{9.5cm}
\begin{center}
  \includegraphics[width=90mm]{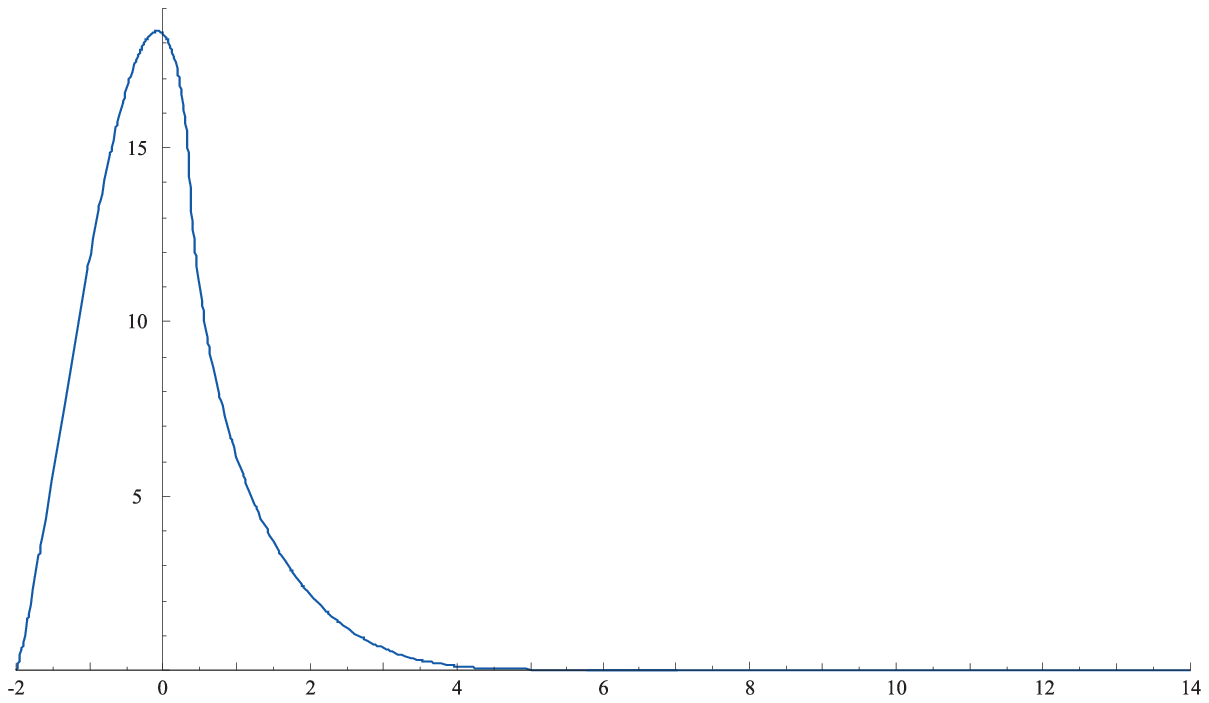}
 \caption{$J_2^{G_2}(y)$} \label{fig-Jy-G2}
\end{center}
\end{minipage}
\end{figure}

We now determine the spectral measure $\nu_{v_N^j}$ over $I_j$.
From (\ref{eqn:measureT2=12C}) and (\ref{eqn:integral_overDj-G2}), with the measure given in Section \ref{sect:measureAinftyG2-T2}, we have that
\begin{equation} \label{eqn:integral_overDj-G2-2}
\frac{1}{192 \pi^4} \int_{\mathbb{T}^2} \psi(\chi_j(\omega_1,\omega_2)) J(\omega_1,\omega_2)^2 \mathrm{d}\omega_1 \, \mathrm{d}\omega_2 = \frac{1}{16 \pi^4} \int_{\mathfrak{D}} \psi(x') |J(x,y)| \mathrm{d}x \, \mathrm{d}y,
\end{equation}
where $\mathfrak{D}$ is as in Section \ref{sect:measureA12inftyG2-D}, and $x'=x,y$ for $j=1,2$ respectively.
Thus the joint spectral measure over $\mathfrak{D}$ is $|J(x,y)| \mathrm{d}x \, \mathrm{d}y/16\pi^4$, which is the reduced Haar measure on $G_2$ \cite[$\S$6.3]{uhlmann/meinel/wipf:2007}.
The measure $\nu_{v_N^1}$ over $I_1$ is obtained by integrating with respect to $y$ in (\ref{eqn:integral_overDj-G2-2}), whilst the measure $\nu_{v_N^2}$ over $I_2$ is obtained by integrating with respect to $x$ in (\ref{eqn:integral_overDj-G2-2}).
More explicitly, using the expressions for the boundaries of $\mathfrak{D}$ given in (\ref{eqn:boundaryD-y-G2-1})-(\ref{eqn:boundaryD-x-G2-3}), the spectral measure $\nu_{v_N^1}$ (over $[-2,7]$) for the graph $\mathcal{A}^{\rho_1}_{\infty}(G_2)$ is $\mathrm{d}\nu_{v_N^1}(x) = J_1^{G_2}(x) \, \mathrm{d}x / 16 \pi^4$, where $J_1^{G_2}(x)$ is given by
$$J_1^{G_2}(x) = \left\{
\begin{array}{cl}
\displaystyle{\int_{-5(x+1)-2(x+2)^{3/2}}^{-5(x+1)+2(x+2)^{3/2}} |J(x,y)| \, \mathrm{d}y} & \textrm{ for } x \in [-2,-1], \\
\displaystyle{\int_{(x^2+2x-7)/4}^{-5(x+1)+2(x+2)^{3/2}} |J(x,y)| \, \mathrm{d}y} & \textrm{ for } x \in [-1,7].
\end{array} \right. $$
The weight $J_1^{G_2}(x)$ is the integral of the square root of a cubic in $y$, and thus can be written in terms of the complete elliptic integrals $K(m)$, $E(m)$ of the first, second kind respectively, where $K(m) = \int_0^{\pi/2} (1-m\sin^2\theta)^{-1/2} \mathrm{d}\theta$ and $E(m) = \int_0^{\pi/2} (1-m\sin^2\theta)^{1/2} \mathrm{d}\theta$. Using \cite[equation 235.14]{byrd/friedman:1971}, $J_1^{G_2}(x)$ is given by
\begin{align*}
\frac{\pi^2}{15} (8(x+2)^{3/2} - x^2-22x-&13)^{1/2} \bigg[ (x^4+236x^3+1662x^2+2876x+1705) \; E(v(x)) \\
&- (8(x+2)^{3/2} + x^2+22x+13)(x^2+22x+13) \; K(v(x)) \bigg],
\end{align*}
for $x \in [-2,-1]$, where $v(x) = 16(x+2)^{3/2}/(8(x+2)^{3/2}-x^2-22x-13)$, whilst for $x \in [-1,7]$, $J_1^{G_2}(x)$ is given by
\begin{align*}
\frac{2\pi^2}{15} (x+2)^{3/4} \bigg[& 2(x^4+236x^3+1662x^2+2876x+1705) \; E(v(x)^{-1}) \\
&- (8(x+2)^{3/2} + x^2+22x+13)(24(x+2)^{3/2} + x^2+22x+13) \; K(v(x)^{-1}) \bigg].
\end{align*}
The weight $J_1^{G_2}(x)$ is illustrated in Figure \ref{fig-Jx-G2}, up to a factor $4\pi^2$.

The spectral measure $\nu_{v_N^2}$ (over $[-2,14]$) for the graph $\mathcal{A}^{\rho_2}_{\infty}(G_2)$ is $\mathrm{d}\nu_{v_N^2}(y) = J_2^{G_2}(y) \, \mathrm{d}y / 16 \pi^4$, where $J_2^{G_2}(y)$ is given by
\begin{eqnarray*}
& \displaystyle{\int_{-1+4p_1(y)+4p_1(y)^2}^{-1+2(y+2)^{1/2}} |J(x,y)| \, \mathrm{d}x} & \textrm{ for } y \in [-2,10/27], \\
& \displaystyle{\int_{-1+4p_1(y)+4p_1(y)^2}^{-1+4p_3(y)+4p_3(y)^2} |J(x,y)| \, \mathrm{d}x + \int_{-1+4p_2(y)+4p_2(y)^2}^{-1+2(y+2)^{1/2}} |J(x,y)| \, \mathrm{d}x} & \textrm{ for } y \in [10/27,5], \\
& \displaystyle{\int_{-1+4p_2(y)+4p_2(y)^2}^{-1+2(y+2)^{1/2}} |J(x,y)| \, \mathrm{d}x} & \textrm{ for } y \in [5,14].
\end{eqnarray*}
A numerical plot of the weight $J_2^{G_2}(y)$ is illustrated in Figure \ref{fig-Jy-G2}, again up to a factor $4\pi^2$.

\section{Joint spectral measures for nimrep graphs associated to $G_2$ modular invariants} \label{sect:measures_nimrepG2}

Suppose $G$ is the nimrep associated to a $G_2$ braided subfactor at some finite level $k$ with vertex set $G_0$.
We define a state $\varphi$ on $\ell^2(G_0)$ by $\varphi( \, \cdot \, ) = \langle \,\cdot \, \Omega, \Omega \rangle$, where $\Omega$ is the basis vector in $\ell^2(G_0)$ corresponding to a distinguished vertex $\ast$. Note that the state $\varphi$ (and thus the spectral measure) depends on the choice of distinguished vertex $\ast$. We choose the distinguished vertex $\ast$ to be the vertex with lowest Perron-Frobenius weight.

Consider the nimrep graph $G_{\lambda}$.
The eigenvalues $\beta_{\lambda}^{(\mu)}$ of $G_{\lambda}$ are given by the ratio $S_{\lambda,\mu}/S_{0,\mu}$, where $\mu$ belongs to the set $\mathrm{Exp}(G)$ of exponents of $G$, $\mathrm{Exp}(G) \subset P^k_{+} = \{ (\lambda_1,\lambda_2) | \, \lambda_1,\lambda_2 \geq 0; \lambda_1 + 2\lambda_2 \leq k \}$ (note that we are now using the Dynkin labels), and the $S$-matrix at level $k$ is given by \cite{kac/peterson:1984, gannon:2001}:
\begin{align*}
S_{\lambda,\mu} & = \frac{-2}{(k+4)\sqrt{3}} \bigg[
\cos(2\xi(2(\hat{\lambda}_1+\hat{\lambda}_2)(\hat{\mu}_1+\hat{\mu}_2) +(\hat{\lambda}_1+\hat{\lambda}_2)\hat{\mu}_2 +\hat{\lambda}_2(\hat{\mu}_1+\hat{\mu}_2) +2\hat{\lambda}_2\hat{\mu}_2)) \\
& \qquad \quad + \cos(2\xi(-(\hat{\lambda}_1+\hat{\lambda}_2)(\hat{\mu}_1+\hat{\mu}_2) -2(\hat{\lambda}_1+\hat{\lambda}_2)\hat{\mu}_2 +\hat{\lambda}_2(\hat{\mu}_1+\hat{\mu}_2) -\hat{\lambda}_2\hat{\mu}_2)) \\
& \qquad \quad + \cos(2\xi(-(\hat{\lambda}_1+\hat{\lambda}_2)(\hat{\mu}_1+\hat{\mu}_2) +(\hat{\lambda}_1+\hat{\lambda}_2)\hat{\mu}_2 -2\hat{\lambda}_2(\hat{\mu}_1+\hat{\mu}_2) -\hat{\lambda}_2\hat{\mu}_2)) \\
& \qquad \quad - \cos(2\xi(-(\hat{\lambda}_1+\hat{\lambda}_2)(\hat{\mu}_1+\hat{\mu}_2) -2(\hat{\lambda}_1+\hat{\lambda}_2)\hat{\mu}_2 -2\hat{\lambda}_2(\hat{\mu}_1+\hat{\mu}_2) -\hat{\lambda}_2\hat{\mu}_2)) \\
& \qquad \quad - \cos(2\xi(2(\hat{\lambda}_1+\hat{\lambda}_2)(\hat{\mu}_1+\hat{\mu}_2) +(\hat{\lambda}_1+\hat{\lambda}_2)\hat{\mu}_2 +\hat{\lambda}_2(\hat{\mu}_1+\hat{\mu}_2) -\hat{\lambda}_2\hat{\mu}_2)) \\
& \qquad \quad - \cos(2\xi(-(\hat{\lambda}_1+\hat{\lambda}_2)(\hat{\mu}_1+\hat{\mu}_2) +(\hat{\lambda}_1+\hat{\lambda}_2)\hat{\mu}_2 +\hat{\lambda}_2(\hat{\mu}_1+\hat{\mu}_2) +2\hat{\lambda}_2\hat{\mu}_2)) \bigg],
\end{align*}
where $\xi=\pi/3(k+4)$, $\lambda=(\lambda_1,\lambda_2)$, $\mu=(\mu_1,\mu_2)$, and $\hat{\lambda}_i=\lambda_i+1$, $\hat{\mu}_i=\mu_i+1$ for $i=1,2$.
Letting $\hat{\mu}_1 = 3(k+4)x_1$, $\hat{\mu}_2 = -(k+4)x_2$, we obtain
\begin{align*}
S_{\lambda,\mu} & = \frac{-2}{(k+4)\sqrt{3}} \bigg[
\cos(2\pi((2\hat{\lambda}_1+3\hat{\lambda}_2)x_1 - (\hat{\lambda}_1+2\hat{\lambda}_2)x_2)) + \cos(2\pi((\hat{\lambda}_1+3\hat{\lambda}_2)x_1 - \hat{\lambda}_2x_2)) \\
& \qquad \qquad + \cos(2\pi(\hat{\lambda}_1x_1 - (\hat{\lambda}_1+\hat{\lambda}_2)x_2)) - \cos(2\pi((\hat{\lambda}_1+3\hat{\lambda}_2)x_1 - (\hat{\lambda}_1+2\hat{\lambda}_2)x_2)) \\
& \qquad \qquad - \cos(2\pi((2\hat{\lambda}_1+3\hat{\lambda}_2)x_1 - (\hat{\lambda}_1+\hat{\lambda}_2)x_2)) - \cos(2\pi(\hat{\lambda}_1x_1 + 2\hat{\lambda}_2x_2)) \bigg],
\end{align*}
which is (up to a scalar factor) nothing but the $S$-function $S_{\lambda + \varrho}(x)$ (see \cite[Section 2.3]{evans/pugh:2012iv} for a discussion on orbit functions).
Then we see that
\begin{equation}\label{eqn:evalues-nimrep-G2}
\beta_{\lambda}^{(\mu)} = \frac{S_{\lambda,\mu}}{S_{0,\mu}} = \frac{S_{\lambda + \varrho}(x)}{S_{\varrho}(x)} = \chi_{\lambda}(t_{\nu}) \in \chi_{\lambda}(\mathbb{T}^2) = I_{\lambda},
\end{equation}
where $t_{\nu}=(e^{2\pi i x_1}, e^{2 \pi i x_2})$, and hence the spectrum $\sigma(G_{\lambda})$ of $G_{\lambda}$ is contained in $I_{\lambda}$. Note that here we are using the Dynkin labels whereas in Section \ref{sect:rep_theoryG2T2} we used labels $(\mu_1,\mu_2) = (\lambda_1+\lambda_2,\lambda_2)$.

Consider now the pair of nimrep graphs $G_{\lambda}$, $G_{\mu}$, which have joint spectrum $\mathfrak{D}_{\lambda,\mu} \subset I_{\lambda} \times I_{\mu}$.
The $m,n^{\mathrm{th}}$ cross moment $\varsigma_{m,n} = \int_{\mathfrak{D}_{\lambda,\mu}} x^m y^n \mathrm{d}\widetilde{\nu}(x,y)$, where $x=x_{\lambda}$, $y=x_{\mu}$, is given by $\langle G_{\lambda}^m G_{\mu}^n \Omega, \Omega \rangle$.
Let $\beta_{\lambda}^{(\sigma)}$ be the eigenvalues of $G_{\lambda}$ with corresponding eigenvectors $\psi^{(\sigma)}$, normalized so that each $\psi^{(\sigma)}$ has norm 1. As the nimreps are a family of commuting matrices they can be simultaneously diagonalised, and thus the eigenvectors of $G_{\lambda}$ are the same for all $\lambda$). Then $G_{\lambda}^m G_{\mu}^n = \mathcal{U} \Lambda_{\lambda}^m \Lambda_{\mu}^n \mathcal{U}^{\ast}$, where $\Lambda_{\lambda}$ is the diagonal matrix with the eigenvalues $\beta_{\lambda}^{(\sigma)}$ on the diagonal, and $\mathcal{U}$ is the unitary matrix whose columns are given by the eigenvectors $\psi^{(\sigma)}$, so that
\begin{equation}\label{eqn:moments-nimrep-G2}
\varsigma_{m,n} \;\; = \;\; \langle \mathcal{U} \Lambda_{\lambda}^m \Lambda_{\mu}^n \mathcal{U}^{\ast} \Omega, \Omega \rangle \;\; = \;\; \langle \Lambda_{\lambda}^m \Lambda_{\mu}^n \mathcal{U}^{\ast} \Omega, \mathcal{U}^{\ast} \Omega \rangle \;\; = \;\; \sum_{\sigma} (\beta_{\lambda}^{(\sigma)})^m (\beta_{\mu}^{(\sigma)})^n |\psi^{(\sigma)}_{\ast}|^2,
\end{equation}
where $\psi^{(\sigma)}_{\ast} = \mathcal{U}^{\ast} \Omega$ is the entry of the eigenvector $\psi^{(\sigma)}$ corresponding to the distinguished vertex $\ast$.
Thus there is a $D_{12}$-invariant measure $\varepsilon$ over $\mathbb{T}^2$ such that
$$\varsigma_{m,n} = \int_{\mathbb{T}^2} \chi_{\lambda}(\omega_1,\omega_2)^m \chi_{\mu}(\omega_1,\omega_2)^n \mathrm{d}\varepsilon(\omega_1,\omega_2),$$
for all $\lambda$, $\mu$.

Note from (\ref{eqn:evalues-nimrep-G2}), (\ref{eqn:moments-nimrep-G2}) that the measure $\varepsilon$ is a discrete measure which has weight $|\psi^{(\nu)}_{\ast}|^2$ at the points $g(t_{\nu}) \in \mathbb{T}^2$ for $g \in D_{12}$, $\nu \in \mathrm{Exp}(G)$, and zero everywhere else. Thus the measure $\varepsilon$ does not depend on the choice of $\lambda$, $\mu$, so that the measure over $\mathbb{T}^2$ is the same for any pair $(G_{\lambda},G_{\mu})$, even though the corresponding measures over $\mathfrak{D}_{\lambda,\mu} \subset \mathbb{R}^2$, and indeed the subsets $\mathfrak{D}_{\lambda,\mu}$ themselves, are different for each such pair.

We will now determine this $D_{12}$-invariant measure $\varepsilon$ over $\mathbb{T}^2$ for all the known $G_2$ modular invariants, where we will focus in particular on the nimrep graphs for the fundamental generators $\rho_j$, $j=1,2$,which have quantum dimensions $[2][7][12]/[4][6]$, $[7][8][15]/[3][4][5]$ respectively, where $[m]$ denotes the quantum integer $[m] = (q^m-q^{-m})/(q-q^{-1})$ for $q=e^{i\pi/3(k+4)}$.
The nimrep graphs $G_{\rho_1}$ were found in \cite{di_francesco:1992}, whilst $G_{\rho_2}$ for the conformal embeddings at levels 3, 4 were found in \cite{coquereaux/rais/tahri:2010}.
The realisation of modular invariants for $G_2$ by braided subfactors is parallel to the realisation of $SU(2)$ and $SU(3)$ modular invariants by $\alpha$-induction for a suitable braided subfactors \cite{ocneanu:2000i, ocneanu:2002, xu:1998, bockenhauer/evans:1999i, bockenhauer/evans:1999ii, bockenhauer/evans/kawahigashi:1999, bockenhauer/evans/kawahigashi:2000}, \cite{ocneanu:2000ii, ocneanu:2002, xu:1998, bockenhauer/evans:1999i, bockenhauer/evans:1999ii, bockenhauer/evans/kawahigashi:1999, bockenhauer/evans:2001, bockenhauer/evans:2002, evans/pugh:2009i, evans/pugh:2009ii} respectively. The realisation of modular invariants for $C_2$ is also under way \cite{evans/pugh:2012iv}.

\subsection{Graphs $\mathcal{A}_k(G_2)$, $k \leq \infty$} \label{sect:measures_AkG2}

The graphs $\mathcal{A}^{\rho_j}_k(G_2)$, $j=1,2$, are associated with the trivial $G_2$ modular invariant at level $k$. They are illustrated in Figures \ref{fig-A_infty(G2)1}, \ref{fig-A_infty(G2)2} respectively, where the set of vertices is now given by $P^k_{++}$. The set of edges is given by the edges between these vertices, except for certain self-loops at the cut-off which are indicated by dashed lines (in Figures \ref{fig-A_infty(G2)1}, \ref{fig-A_infty(G2)2} the dashed lines indicate the edges to be removed when $k=6$).
The eigenvalues $\beta^{j,(\lambda)} := \beta_{\rho_j}^{(\lambda)}$ of $\mathcal{A}^{\rho_j}_k(G_2)$, $j=1,2$, are given by the ratio $S_{\rho_j,\lambda}/S_{0,\lambda}$ with corresponding eigenvectors $\psi^{\lambda}_{\mu} = S_{\lambda,\mu}$, where $\mu \in \mathrm{Exp}(\mathcal{A}_k(G_2)) = P^k_{++}$. Then with $\mu = \ast = (0,0)$, we obtain
\begin{eqnarray}
\psi^{\lambda}_{\ast} & = & \frac{-2}{(k+4)\sqrt{3}} \bigg[ \cos(2\xi(5\hat{\lambda}_1+9\hat{\lambda}_2)) + \cos(2\xi(\hat{\lambda}_1+6\hat{\lambda}_2)) + \cos(2\xi(4\hat{\lambda}_1+3\hat{\lambda}_2)) \nonumber \\
&& \qquad \qquad \quad - \cos(2\xi(4\hat{\lambda}_1+9\hat{\lambda}_2)) - \cos(2\xi(\hat{\lambda}_1-3\hat{\lambda}_2)) - \cos(2\xi(5\hat{\lambda}_1+6\hat{\lambda}_2)) \bigg] \qquad \label{eqn:PF-evector1-G2}
\end{eqnarray}
and hence we see that
\begin{equation} \label{eqn:PF-evector=J-G2}
\psi^{\lambda}_{\ast} = \frac{-1}{4 \sqrt{3} (k+4) \pi^2} J \left( (\hat{\lambda}_1 + 3\hat{\lambda}_2)/3(k+4), -\hat{\lambda}_1/3(k+4) \right),
\end{equation}
where $J(\theta_1,\theta_2)$ is given by (\ref{eqn:J[theta]-G2}) and in (\ref{eqn:PF-evector=J-G2}) we have
\begin{equation} \label{eqn:theta-lambda_G2}
\theta_1 = (\hat{\lambda}_1 + 3\hat{\lambda}_2)/3(k+4), \qquad \theta_2 = -\hat{\lambda}_1/3(k+4),
\end{equation}
so that $\hat{\lambda}_1 = -3(k+4)\theta_2$ and $\hat{\lambda}_2 = (k+4)(\theta_1+\theta_2)$.

As a consequence of the identification (\ref{eqn:PF-evector=J-G2}) between the Perron-Frobenius eigenvector and the Jacobian we can obtain another expression for the Jacobian $J$. Recall that the Perron-Frobenius eigenvector for $\mathcal{A}_k(G_2)$ can also be written in the Kac-Weyl factorized form \cite{kac/peterson:1984}:
\begin{equation} \label{eqn:PF-evector2-G2}
\phi^{\ast}_{\lambda} = \frac{ \sin(\hat{\lambda}_1\xi) \sin(3\hat{\lambda}_2\xi) \sin((\hat{\lambda}_1+3\hat{\lambda}_2)\xi) \sin((2\hat{\lambda}_1+3\hat{\lambda}_2)\xi) \sin((3\hat{\lambda}_1+3\hat{\lambda}_2)\xi) \sin((3\hat{\lambda}_1+6\hat{\lambda}_2)\xi) }{ \sin(\xi) \sin(3\xi) \sin(4\xi) \sin(5\xi) \sin(6\xi) \sin(9\xi) }.
\end{equation}
Now $\phi^{\ast}_{\ast} = 1$ whilst $\psi_{\ast}^{\ast} = -2[\cos(28\xi)+2\cos(14\xi)-\cos(26\xi)-\cos(22\xi)-\cos(4\xi)]/(k+4)\sqrt{3} = 64\sin(\xi) \sin(3\xi) \sin(4\xi) \sin(5\xi) \sin(6\xi) \sin(9\xi)/(k+4)\sqrt{3}$. Thus we see that $(k+4)\sqrt{3} \psi_{\ast} = 64\sin(\xi) \sin(3\xi) \sin(4\xi) \sin(5\xi) \sin(6\xi) \sin(9\xi)\phi^{\ast}$. Then from (\ref{eqn:PF-evector=J-G2}) we have
\begin{eqnarray*}
\lefteqn{ J(\theta_1,\theta_2) \;\; = \;\; - 4 (k+4) \sqrt{3} \pi^2 \; \psi_{\ast}^{(-3(k+4)(\theta_2-1),(k+4)(\theta_1+\theta_2-1))} } \\
& = & - 256 \pi^2 \sin(\xi) \sin(3\xi) \sin(4\xi) \sin(5\xi) \sin(6\xi) \sin(9\xi) \; \phi^{\ast}_{(-3(k+4)(\theta_2-1),(k+4)(\theta_1+\theta_2-1))} \\
& = & 256 \pi^2 \sin(\theta_1\pi) \sin(\theta_2\pi) \sin((\theta_1+\theta_2)\pi) \sin((\theta_1-\theta_2)\pi) \sin((2\theta_1-\theta_2)\pi) \sin((\theta_1-2\theta_2)\pi),
\end{eqnarray*}
so that the Jacobian $J(\theta_1,\theta_2)$ can also be written as a product of sine functions.

The eigenvalues $\beta^{j,(\lambda)} = S_{\rho_j,\lambda}/S_{0,\lambda}$ are given by
\begin{eqnarray*}
\beta^{1,(\lambda)} & = & 1 + 2\cos(2\xi\hat{\lambda}_1) + 2\cos(2\xi(\hat{\lambda}_1+3\hat{\lambda}_2)) + 2\cos(2\xi(2\hat{\lambda}_1+3\hat{\lambda}_2)) \;\; = \;\; \chi_1(\omega_1,\omega_2), \\
\beta^{2,(\lambda)} & = & 2 + 2\cos(2\xi\hat{\lambda}_1) + 2\cos(2\xi(\hat{\lambda}_1+3\hat{\lambda}_2)) + 2\cos(2\xi(2\hat{\lambda}_1+3\hat{\lambda}_2)) \\
&& + 2\cos(6\xi\hat{\lambda}_2) + 2\cos(6\xi(\hat{\lambda}_1+\hat{\lambda}_2)) + 2\cos(6\xi(\hat{\lambda}_1+2\hat{\lambda}_2)) \qquad \;\; = \;\; \chi_2(\omega_1,\omega_2),
\end{eqnarray*}
where $\omega_j=\exp^{2\pi i \theta_j}$, $j=1,2$ are related to $\lambda$ as in (\ref{eqn:theta-lambda_G2}).

We now compute the spectral measure for $\mathcal{A}^{\rho_j}_k(G_2)$.
Now summing over all $(\lambda_1,\lambda_2) \in \mathrm{Exp}(\mathcal{A}_k(G_2))$ corresponds to summing over all $(\theta_1, \theta_2) \in \{ ((\hat{\lambda}_1+3\hat{\lambda}_2)/3(k+4),-\hat{\lambda}_1/3(k+4)) | \; \hat{\lambda}_1,\hat{\lambda}_2 \geq 1, \hat{\lambda}_1+2\hat{\lambda}_2 \leq k+3 \}$, or equivalently, since such points satisfy $\theta_1 + \theta_2 \equiv 0 \textrm{ mod } 3$, to summing over all $(\theta_1, \theta_2) \in F_k' = \{ (q_1/3(k+4), q_2/3(k+4)) | \; q_1, q_2 = 0,1,\ldots,3k+11; \, q_1 + q_2 \equiv 0 \textrm{ mod } 3 \}$ such that
\begin{align*}
\theta_2 & = -\hat{\lambda}_1/3(k+4) \leq -1/3(k+4), \qquad
\theta_1 + \theta_2 = \hat{\lambda}_2/(k+4) \geq 1/(k+4) \\
2\theta_1 - \theta_2 & = (\hat{\lambda}_1 + 2\hat{\lambda}_2)/(k+4) \leq (k+3)/(k+4) = 1 - 1/(k+4).
\end{align*}
Since $\theta_j$ and $\theta_j + n$ give the same points in $\mathbb{T}^2$ for $n \in \mathbb{Z}$, the last three conditions are equivalent to
$$\theta_2 \leq 1-1/3(k+4), \qquad \theta_1 + \theta_2 \geq 1+1/(k+4), \qquad 2\theta_1 - \theta_2 \leq - 1/(k+4).$$

Denote by $F_k$ the set of all $(\omega_1,\omega_2) \in \mathbb{T}^2$ such that $(\theta_1,\theta_2) \in F_k'$ satisfies these conditions. Then from (\ref{eqn:moments-nimrep-G2}) and (\ref{eqn:PF-evector=J-G2}) we obtain
\begin{eqnarray}
\varsigma_{m,n} & = & \frac{1}{48(k+4)^2 \pi^4} \sum_{\lambda \in \mathrm{Exp}} (\beta^{1,(\lambda)})^m (\beta^{2,(\lambda)})^n J \left( (\lambda_1 + 3\lambda_2)/3(k+4), -\lambda_1/3(k+4) \right)^2 \nonumber \\
& = & \frac{1}{48(k+4)^2 \pi^4} \sum_{(\omega_1,\omega_2) \in F_k} (\chi_1(\omega_1,\omega_2))^m (\chi_2(\omega_1,\omega_2))^n J(\omega_1,\omega_2)^2 \label{eqn:sum_for_measure-A(G2)}
\end{eqnarray}

\begin{figure}[tb]
\begin{center}
  \includegraphics[width=55mm]{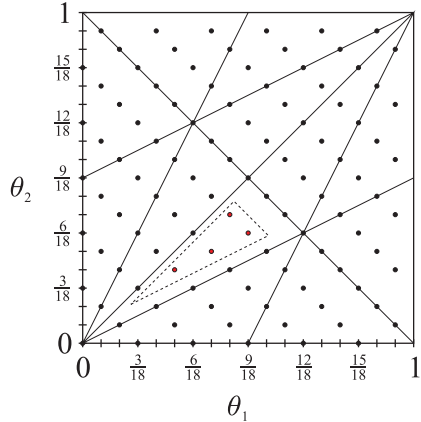}\\
 \caption{The points $(\theta_1,\theta_2)$ such that $(e^{2 \pi i \theta_1}, e^{2 \pi i \theta_2}) \in F_2^W$.} \label{fig:C2W-G2}
\end{center}
\end{figure}

If we let $F$ be the limit of $F_k$ as $k \rightarrow \infty$, then $F$ is a fundamental domain of $\mathbb{T}^2$ under the action of the group $D_{12}$, illustrated in Figure \ref{fig:fund_domain-G2inT2}. Since $J=0$ along the boundary of $F$, which is mapped to the boundary of $\mathfrak{D}$ under $\Psi:\mathbb{T}^2 \rightarrow \mathfrak{D}$, we can include points on the boundary of $F$ in the summation in (\ref{eqn:sum_for_measure-A(G2)}). Since $J^2$ is invariant under the action of $D_{12}$, we have
\begin{eqnarray}
\varsigma_{m,n} & = & \frac{1}{12} \frac{1}{48(k+4)^2 \pi^4} \sum_{(\omega_1,\omega_2) \in F_k^W} (\chi_1(\omega_1,\omega_2))^m (\chi_2(\omega_1,\omega_2))^n J(\omega_1,\omega_2)^2 \label{eqn:sum_for_measure-A(G2)2}
\end{eqnarray}
where
\begin{equation} \label{def:Dl}
F_k^W = \{ (e^{2 \pi i q_1/3(k+4)}, e^{2 \pi i q_2/3(k+4)}) \in \mathbb{T}^2 | \; q_1,q_2 = 0, 1, \ldots, 3k+11; \, q_1 + q_2 \equiv 0 \textrm{ mod } 3 \}
\end{equation}
is the image of $F_k$ under the action of the Weyl group $W=D_{12}$ and some additional points which lie on the boundaries of the fundamental domains (i.e. where $J=0$). We illustrate the points $(\theta_1,\theta_2)$ such that $(e^{2 \pi i \theta_1}, e^{2 \pi i \theta_2}) \in F_2^W$ in Figure \ref{fig:C2W-G2}. The points in the interior of the fundamental domain $F$, those enclosed by the dashed line, correspond to the vertices of the graph $\mathcal{A}_2(G_2)$.

Note that $F_k^W = D_{k+4}$ in the notation of \cite[$\S$7.1]{evans/pugh:2009v}, and that $|F_k^W| = 3(k+4)^2$.
Thus from (\ref{eqn:sum_for_measure-A(G2)2}), we obtain (c.f. \cite[Theorem 4]{evans/pugh:2009v}):

\begin{Thm}
The joint spectral measure of $\mathcal{A}^{\rho_j}_k(G_2)$, $j=1,2$, (over $\mathbb{T}^2$) is given by
\begin{equation}
\mathrm{d}\varepsilon(\omega_1,\omega_2) = \frac{1}{192 \pi^4} J(\omega_1,\omega_2)^2 \mathrm{d}^{(k+4)}(\omega_1,\omega_2),
\end{equation}
where $\mathrm{d}^{(k+4)}$ is the uniform measure over $F_k^W$.
\end{Thm}

In fact, the spectral measure over $\mathbb{T}^2$ for the nimrep graph $G_{\lambda}$ for any $\lambda \in {}_N \mathcal{X}_N$ (where $G_{\rho_j} = \mathcal{A}^{\rho_j}_k(G_2)$) is given by the above measure.

We can now easily deduce the spectral measure (over $\mathbb{T}^2$) for $\mathcal{A}_{\infty}(G_2)$ claimed in Section \ref{sect:measureAinftyG2-T2}. Letting $k \rightarrow \infty$, the measure $\mathrm{d}^{(k+4)}(\omega_1,\omega_2)$ becomes the uniform Lebesgue measure $\mathrm{d}\omega_1 \, \mathrm{d}\omega_2$ on $\mathbb{T}^2$. Thus:

\begin{Thm}
The joint spectral measure of the infinite graph $\mathcal{A}_{\infty}(G_2)$ (over $\mathbb{T}^2$) is given by
\begin{equation}
\mathrm{d}\varepsilon(\omega_1,\omega_2) = \frac{1}{192 \pi^4} J(\omega_1,\omega_2)^2 \mathrm{d}\omega_1 \, \mathrm{d}\omega_2,
\end{equation}
where $\mathrm{d}\omega$ is the uniform Lebesgue measure over $\mathbb{T}$.
\end{Thm}

Then the spectral measure for $\mathcal{A}^{\rho_j}_k(G_2)$ over $\mathfrak{D}$ or $I_j$, $j=1,2$, has the same weights as the spectral measure for the infinite graph $\mathcal{A}^{\rho_j}_{\infty}(G_2)$ given in Section \ref{sect:measureAinftyG2-R}, but the measure here is a discrete measure.

\renewcommand{\arraystretch}{1.5}

\subsection{Exceptional Graph $\mathcal{E}_3(G_2)$: $(G_2)_3 \rightarrow (E_6)_1$} \label{sect:measures_E3G2}

\begin{figure}[tb]
\begin{center}
  \includegraphics[width=25mm]{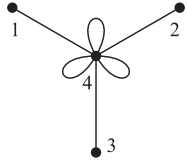} \\
 \caption{Exceptional Graph $\mathcal{E}_3(G_2)$} \label{fig-Graph_E3_G2}
\end{center}
\end{figure}

The graph $\mathcal{E}_3(G_2) := \mathcal{E}_3^{\rho_1}(G_2) = \mathcal{E}_3^{\rho_2}(G_2)$, illustrated in Figure \ref{fig-Graph_E3_G2}, is the nimrep graph associated to the conformal embedding $(G_2)_3 \rightarrow (E_6)_1$, and is one of two nimrep graphs associated to the modular invariant
$$Z_{\mathcal{E}_3} = |\chi_{(0,0)}+\chi_{(1,1)}|^2 + 2|\chi_{(2,0)}|^2$$
which is at level 3 and has exponents $\mathrm{Exp}(\mathcal{E}_3(G_2)) = \{ (0,0), (1,1), \textrm{ and } (2,0) \textrm{ twice } \}$. The other nimrep graph associated to this modular invariant is $\mathcal{E}_3^M(G_2)$ considered in the next section.

\begin{figure}[tb]
\begin{center}
  \includegraphics[width=75mm]{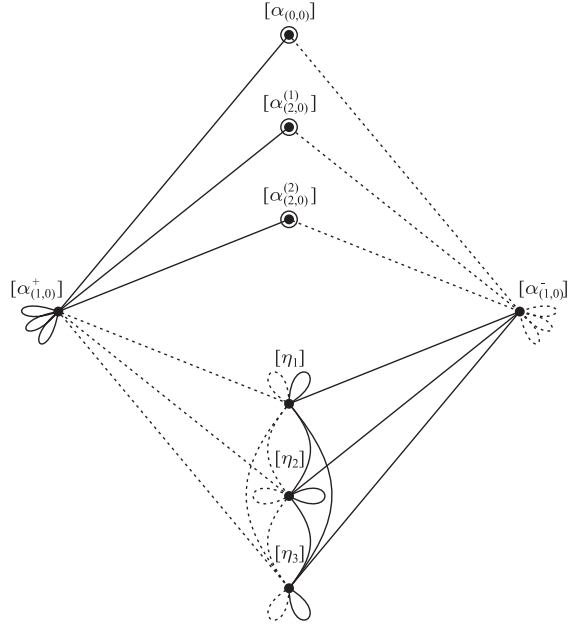} \\
 \caption{$\mathcal{E}_3(G_2)$: Multiplication by $[\alpha_{(1,0)}^+]$ (solid lines) and $[\alpha_{(1,0)}^-]$ (dashed lines)} \label{Fig-full_E3_G2}
\end{center}
\end{figure}

Following \cite[$\S$6]{bockenhauer/evans:1999ii} we can compute the principal graph and dual principal graph of the inclusion $(G_2)_3 \rightarrow (E_6)_1$. The chiral induced sector bases ${}_M \mathcal{X}_M^{\pm} \subset \mbox{Sect}(M)$ and full induced sector basis ${}_M \mathcal{X}_M \subset \mbox{Sect}(M)$, the sector bases given by all irreducible subsectors of $[\alpha_{\lambda}^{\pm}]$ and $[\alpha_{\lambda}^+ \circ \alpha_{\lambda'}^-]$ respectively, for $\lambda, \lambda' \in {}_N \mathcal{X}_N$, are given by
\begin{align*}
{}_M \mathcal{X}_M^{\pm} &= \{ [\alpha_{(0,0)}], [\alpha_{(1,0)}^{\pm}], [\alpha_{(2,0)}^{(1)}], [\alpha_{(2,0)}^{(2)}] \},\\
{}_M \mathcal{X}_M &= \{ [\alpha_{(0,0)}], [\alpha_{(1,0)}^+], [\alpha_{(1,0)}^-], [\alpha_{(2,0)}^{(1)}], [\alpha_{(2,0)}^{(2)}], [\eta_1], [\eta_2], [\eta_3] \},
\end{align*}
where $[\alpha_{(2,0)}^{\pm}] = [\alpha_{(1,0)}^{\pm}] \oplus [\alpha_{(2,0)}^{(1)}] \oplus [\alpha_{(2,0)}^{(2)}]$, $[\alpha_{(1,0)}^+ \circ \alpha_{(1,0)}^-] = [\eta_1] \oplus [\eta_2] \oplus [\eta_3]$, and $\alpha_{(i,j)} \equiv \alpha_{\lambda_{(i,j)}}$.
The fusion graphs of $[\alpha_{(1,0)}^+]$ (solid lines) and $[\alpha_{(1,0)}^-]$ (dashed lines) are given in Figure \ref{Fig-full_E3_G2}, see also \cite[Figure 17(a)]{coquereaux/rais/tahri:2010}. The marked vertices corresponding to sectors in ${}_M \mathcal{X}_M^0 = {}_M \mathcal{X}_M^+ \cap {}_M \mathcal{X}_M^-$ have been circled.
Note that multiplication by $[\alpha_{(1,0)}^+]$ (or $[\alpha_{(1,0)}^-]$) does not give two copies of the nimrep graph $\mathcal{E}_3(G_2)$ as one might expect, but rather one copy each of $\mathcal{E}_3(G_2)$ and $\mathcal{E}_3^M(G_2)$. This is similar to the situation for the $SU(3)$ conformal embedding $SU(3)_9 \rightarrow (E_6)_1$ \cite[$\S$5.2]{evans:2002}.

\begin{figure}[tb]
\begin{minipage}[t]{7.5cm}
\begin{center}
  \includegraphics[width=62.5mm]{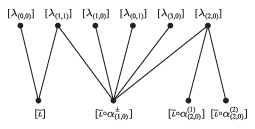} \\
 \caption{$\mathcal{E}_3(G_2)$: Principal graph of $(G_2)_3 \rightarrow (E_6)_1$} \label{Fig-GHJ_Graph_E3_G2}
\end{center}
\end{minipage}
\hfill
\begin{minipage}[t]{7.5cm}
\begin{center}
  \includegraphics[width=75mm]{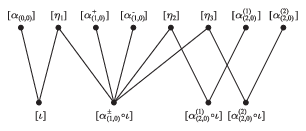} \\
 \caption{$\mathcal{E}_3(G_2)$: Dual principal graph of $(G_2)_3 \rightarrow (E_6)_1$} \label{Fig-GHJ_Graph_E3_G2_dual}
\end{center}
\end{minipage}
\end{figure}

Let $\iota:N \hookrightarrow M$ denote the injection map $\iota(n)=n \in M$, $n \in N$ and $\overline{\iota}$ its conjugate. The dual canonical endomorphism $\theta = \overline{\iota} \circ \iota$ for the conformal embedding can be read from the vacuum block of the modular invariant: $[\theta] = [\lambda_{(0,0)}] \oplus [\lambda_{(1,1)}]$.
By \cite[Corollary 3.19]{bockenhauer/evans:1999ii} and the fact that $\langle \gamma, \gamma \rangle_M = \langle \theta, \theta \rangle_N = 2$, the canonical endomorphism $\gamma = \iota \circ \overline{\iota}$ is given by
\begin{equation}
[\gamma] = [\alpha_{(0,0)}] \oplus [\eta_1].
\end{equation}
Then by \cite[Theorem 4.2]{bockenhauer/evans:1999ii}, the principal graph of the inclusion $(G_2)_3 \rightarrow (E_6)_1$ of index $\frac{1}{2}(7+\sqrt{21}) \approx 5.79$ is given by the connected component of $[\lambda_{(0,0)}] \in {}_N \mathcal{X}_N$ of the induction-restriction graph, and the dual principal graph is given by the connected component of $[\alpha_{(0,0)}] \in {}_M \mathcal{X}_M$ of the $\gamma$-multiplication graph. The principal graph and dual principal graph are illustrated in Figures \ref{Fig-GHJ_Graph_E3_G2} and \ref{Fig-GHJ_Graph_E3_G2_dual} respectively. These principal graphs are sometimes referred to as ``Haagerup with legs'' \cite[$\S$4.2.4]{jones/morrison/snyder:2013}. The principal graph in Figure \ref{Fig-GHJ_Graph_E3_G2} appears as the intertwiner for the quantum subgroup $\mathcal{E}_3(G_2)$ in \cite[$\S$4.4]{coquereaux/rais/tahri:2010}.

\begin{figure}[tb]
\begin{center}
  \includegraphics[width=45mm]{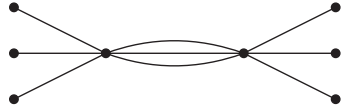} \\
 \caption{$\mathcal{E}_3(G_2)$: Principal graph of $\alpha_{(1,0)}^{\pm}(M) \subset M$} \label{Fig-principal_graph_E3_G2}
\end{center}
\end{figure}

One can also construct a subfactor $\alpha_{(1,0)}^{\pm}(M) \subset M$ with index $(\frac{1}{2}(3+\sqrt{21}))^2 = \frac{3}{2}(5+\sqrt{21}) \approx 14.37$, where $M$ is a type III factor. This subfactor has already appeared in \cite{izumi:2001, xu:2001} (see also the Appendix in \cite{calegari/morrison/snyder:2010}, \cite{jones/morrison/snyder:2013} and \cite{evans/gannon:2012}). The chiral systems ${}_M \mathcal{X}_M^{\pm}$ are near group $C^{\ast}$-category of type $G+3$, where $G = \mathbb{Z}_3$, generated by a self-conjugate irreducible endomorphism $\rho$ of $M$ and an outer action $\alpha$ of $G$ on $M$, such that $[\alpha_i][\rho]=[\rho]=[\rho][\alpha_i]$ and $[\rho^2]=\bigoplus_{i=1}^3 [\alpha_i] \oplus 3[\rho]$, for $i \in \mathbb{Z}_3$. Here $\alpha_{(1,0)}^{\pm} = \rho$, $\alpha_{(0,0)} = \alpha_0$ and $\alpha_{(2,0)}^{(j)} = \alpha_j$. The index $d_{\rho}^2$ of $\rho(M) \subset M$ is thus $d_{\rho}^2=3d_{\rho}+3$, where $d_{\sigma}$ is the statistical dimension of $\sigma$.
Its principal graph is illustrated in Figure \ref{Fig-principal_graph_E3_G2} and is the bipartite unfolded version of the graph $\mathcal{E}_3(G_2)$. The dual principal graph is isomorphic to the principal graph as abstract graphs \cite[Corollary 3.7]{xu:1998}.

We now determine the joint spectral measure of $\mathcal{E}_3^{\rho_1}(G_2)$, $\mathcal{E}_3^{\rho_2}(G_2)$. With $\theta_1, \theta_2$ as in (\ref{eqn:theta-lambda_G2}) for $\lambda = (\lambda_1,\lambda_2) \in \mathrm{Exp}(\mathcal{E}_3(G_2))$, we have the following values:
\begin{center}
\begin{tabular}{|c|c|c|c|} \hline
$\lambda \in \mathrm{Exp}$ & $(\theta_1,\theta_2) \in [0,1]^2$ & $|\psi^{\lambda}_{\ast}|^2$ & $\frac{1}{8\pi^2}|J(\theta_1,\theta_2)|$ \\
\hline $(0,0)$ & $\left(\frac{4}{21},\frac{20}{21}\right)$ & $\frac{7-\sqrt{21}}{42}$ & $\frac{7-\sqrt{21}}{4}$ \\
\hline $(1,1)$ & $\left(\frac{8}{21},\frac{19}{21}\right)$ & $\frac{7+\sqrt{21}}{42}$ & $\frac{7+\sqrt{21}}{4}$ \\
\hline $(2,0)$ & $\left(\frac{2}{7},\frac{6}{7}\right)$ & $\frac{1}{2}$ & $\frac{7}{2}$ \\
\hline
\end{tabular}
\end{center}
where the eigenvectors $\psi^{\lambda}$ have been normalized so that $||\psi^{\lambda}|| = 1$, and for the exponent $(2,0)$ which has multiplicity two, the value listed in the table for $|\psi^{(2,0)}_{\ast}|^2$ is $|\psi^{(2,0)_1}_{\ast}|^2 + |\psi^{(2,0)_2}_{\ast}|^2$.
Note that
\begin{equation} \label{eqn:psi=J+zeta}
|\psi^{\lambda}_{\ast}|^2 = \frac{2}{21} \left( \frac{1}{8\pi^2} |J| \right) + \zeta_{\lambda}
\end{equation}
where $\zeta_{\lambda} = 0$ for $\lambda \in \{ (0,0), (1,1) \}$ and $\zeta_{(2,0)} = 1/6$.

\begin{figure}[tb]
\begin{center}
  \includegraphics[width=45mm]{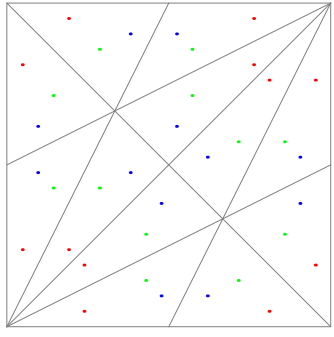}\\
 \caption{The orbit of the points $(\theta_1,\theta_2)$ for $\lambda \in \mathrm{Exp}(\mathcal{E}_3(G_2))$.} \label{Fig-E3pointsG2}
\end{center}
\end{figure}

The orbit under $D_{12}$ of the points $(\theta_1,\theta_2) \in \left\{ \left(\frac{4}{21},\frac{20}{21}\right)\textcolor{red}{\cdot}, \left(\frac{8}{21},\frac{19}{21}\right)\textcolor{blue}{\cdot}, \left(\frac{2}{7},\frac{6}{7}\right)\textcolor{green}{\cdot} \right\}$ are illustrated in Figure \ref{Fig-E3pointsG2}.
These points give the measure $\mathrm{d}^{(21/4,1/21)}$, whose support has cardinality 36. Note that when taking the orbit under $D_{12}$, the associated weight in (\ref{eqn:psi=J+zeta}) is now counted 12 times, thus we must divide (\ref{eqn:psi=J+zeta}) by 12.
Thus the measure for $\mathcal{E}_3(G_2)$ is
$$\mathrm{d}\varepsilon = 36 \, \frac{1}{12} \, \frac{2}{21} \, \frac{1}{8\pi^2} |J| \, \mathrm{d}^{(21/4,1/21)} + \frac{\zeta_{(2,0)}}{12} \sum_{g \in D_{12}} \delta_{g(e^{4\pi i/7},e^{6\pi i/7})},$$
where $\delta_x$ is the Dirac measure at the point $x$. Then we have obtained the following result:

\begin{Thm}
The joint spectral measure of $\mathcal{E}_3^{\rho_1}(G_2)$, $\mathcal{E}_3^{\rho_2}(G_2)$ (over $\mathbb{T}^2$) is
\begin{equation}
\mathrm{d}\varepsilon = \frac{1}{28\pi^2} |J| \, \mathrm{d}^{(21/4,1/21)} + \frac{1}{72} \sum_{g \in D_{12}} \delta_{g(e^{4\pi i/7},e^{6\pi i/7})},
\end{equation}
where $\mathrm{d}^{(n,k)}$ is as in Definition \ref{def:4measures} and $\delta_x$ is the Dirac measure at the point $x$.
\end{Thm}

\subsection{Exceptional Graph $\mathcal{E}_3^M(G_2)$: $(G_2)_3 \rightarrow (E_6)_1 \rtimes \mathbb{Z}_3$} \label{sect:measures_E3MG2}

\begin{figure}[tb]
\begin{center}
  \includegraphics[width=30mm]{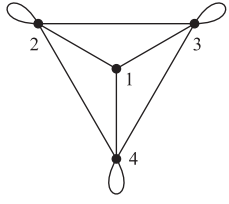} \\
 \caption{Exceptional Graph $\mathcal{E}_3^M(G_2)$} \label{fig-Graph_E3M_G2}
\end{center}
\end{figure}

The graph $\mathcal{E}_3^M(G_2):= \mathcal{E}_3^{M,\rho_1}(G_2) = \mathcal{E}_3^{M,\rho_2}(G_2)$, illustrated in Figure \ref{fig-Graph_E3M_G2}, is the nimrep graph for the type II inclusion $(G_2)_3 \rightarrow (E_6)_1 \rtimes_{\tau} \mathbb{Z}_3$ with index $\frac{3}{2}(7+\sqrt{21}) \approx 17.37$, where $\tau = \alpha_{(2,0)}^{(1)}$ is a non-trivial simple current of order 3 in the ambichiral system ${}_M \mathcal{X}_M^0$, see Figure \ref{Fig-full_E3_G2}. For such an orbifold inclusion to exist, one needs an automorphism $\tau_0$ such that $[\tau_0]=[\tau]$ and $\tau_0^3=\mbox{id}$ \cite[$\S3$]{bockenhauer/evans:1999i}, which exists precisely when the statistics phase $\omega_{\tau}$ of $\tau$ satisfies $\omega_{\tau}^3 = 1$ \cite[Lemma 4.4]{rehren:1990}. By \cite[Lemma 6.1]{bockenhauer/evans:2000}, if $[\tau]$ is a subsector of $[\alpha_{\lambda}^+]$ and $[\alpha_{\mu}^-]$ for some $\lambda, \mu \in {}_N \mathcal{X}_N$, then $\omega_{\tau} = \omega_{\lambda} = \omega_{\mu}$, and hence it is sufficient to check that $\omega_{\lambda}$ and $\omega_{\mu}$ satisfy $\omega^3 = 1$. From Section \ref{sect:measures_E3G2}, $[\tau]$ ($=[\alpha_{(2,0)}^{(1)}]$) is a subsector of $[\alpha_{(2,0)}^{\pm}]$. Now $\omega_{(2,0)} = e^{4\pi i/3}$ \cite[$\S4.4$]{coquereaux/rais/tahri:2010}, which satisfies $\omega_{(2,0)}^3 = 1$, as required.

\begin{figure}[tb]
\begin{center}
  \includegraphics[width=60mm]{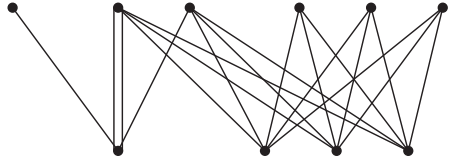} \\
 \caption{$\mathcal{E}_3^M(G_2)$: Principal graph of $(G_2)_3 \rightarrow (E_6)_1 \rtimes \mathbb{Z}_3$} \label{Fig-GHJ_Graph_E3M_G2}
\end{center}
\end{figure}

The principal graph for this inclusion is illustrated in Figure \ref{Fig-GHJ_Graph_E3M_G2}.
This will be discussed in a future publication using a generalised Goodman-de la Harpe-Jones construction analogous to that for the $D_{\mathrm{odd}}$ and $E_7$ modular invariants for $SU(2)$ \cite[$\S5.2, 5.3$]{bockenhauer/evans/kawahigashi:2000} and the type II inclusions for $SU(3)$ \cite[$\S5$]{evans/pugh:2009ii}.
It is not clear what the dual principal graph is in this case.

The associated modular invariant is again $Z_{\mathcal{E}_3}$ and the graph $\mathcal{E}_3^M(G_2)$ is isospectral to $\mathcal{E}_3(G_2)$. In fact, $\mathcal{E}_3^M(G_2)$ is obtained from $\mathcal{E}_3(G_2)$ by a $\mathbb{Z}_3$-orbifold procedure.
Then with $\theta_1, \theta_2$ as in (\ref{eqn:theta-lambda_G2}) for $\lambda = (\lambda_1,\lambda_2) \in \mathrm{Exp}(\mathcal{E}_3^M(G_2)) = \mathrm{Exp}(\mathcal{E}_3(G_2))$, we have:
\begin{center}
\begin{tabular}{|c|c|c|c|} \hline
$\lambda \in \mathrm{Exp}$ & $(\theta_1,\theta_2) \in [0,1]^2$ & $|\psi^{\lambda}_{\ast}|^2$ & $\frac{1}{64\pi^4}J(\theta_1,\theta_2)^2$ \\
\hline $(0,0)$ & $\left(\frac{4}{21},\frac{20}{21}\right)$ & $\frac{7-\sqrt{21}}{14}$ & $\frac{7(5-\sqrt{21})}{8}$ \\
\hline $(1,1)$ & $\left(\frac{8}{21},\frac{19}{21}\right)$ & $\frac{7+\sqrt{21}}{14}$ & $\frac{7(5+\sqrt{21})}{8}$ \\
\hline $(2,0)$ & $\left(\frac{2}{7},\frac{6}{7}\right)$ & 0 & $\frac{49}{4}$ \\
\hline
\end{tabular}
\end{center}
where the eigenvectors $\psi^{\lambda}$ have been normalized so that $||\psi^{\lambda}|| = 1$.
In this case $\zeta_{(2,0)} = -1$ in (\ref{eqn:psi=J+zeta}).
Thus we have the following result:

\begin{Thm}
The joint spectral measure of $\mathcal{E}_3^{M,\rho_1}(G_2)$, $\mathcal{E}_3^{M,\rho_2}(G_2)$ (over $\mathbb{T}^2$) is
\begin{equation}
\mathrm{d}\varepsilon = \frac{3}{28\pi^2} |J| \, \mathrm{d}^{(21/4,1/21)} - \frac{1}{12} \sum_{g \in D_{12}} \delta_{g(e^{4\pi i/7},e^{6\pi i/7})},
\end{equation}
where $\mathrm{d}^{(n,k)}$ are as in Definition \ref{def:4measures} and $\delta_x$ is the Dirac measure at the point $x$.
\end{Thm}

\subsection{Exceptional Graph $\mathcal{E}_4(G_2)$: $(G_2)_4 \rightarrow (D_7)_1$} \label{sect:measures_E4G2}

\begin{figure}[tb]
\begin{minipage}[t]{5.5cm}
\begin{center}
  \includegraphics[width=35mm]{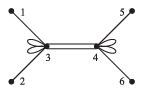} \\
 \caption{Graph $\mathcal{E}_4^{\rho_1}(G_2)$} \label{fig-Graph_E4_G2-1}
\end{center}
\end{minipage}
\hfill
\begin{minipage}[t]{8cm}
\begin{center}
  \includegraphics[width=80mm]{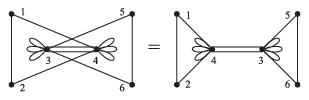} \\
 \caption{Graph $\mathcal{E}_4^{\rho_2}(G_2)$} \label{fig-Graph_E4_G2-2}
\end{center}
\end{minipage}
\end{figure}

The graphs $\mathcal{E}_4^{\rho_j}(G_2)$, illustrated in Figures \ref{fig-Graph_E4_G2-1} and \ref{fig-Graph_E4_G2-2}, are the nimrep graphs associated with the conformal embedding $(G_2)_4 \rightarrow (D_7)_1 = (\mathrm{Spin}(14))_1$ and are one of two families of graphs associated to the modular invariant
$$Z_{\mathcal{E}_4} = |\chi_{(0,0)}+\chi_{(3,0)}|^2 + |\chi_{(0,1)}+\chi_{(4,0)}|^2 + 2|\chi_{(1,1)}|^2$$
at level 4 with exponents $\mathrm{Exp}(\mathcal{E}_4(G_2)) = \{ (0,0), (3,0), (0,1), (4,0) \textrm{ and } (1,1) \textrm{ twice } \}$.

\begin{figure}[tb]
\begin{center}
  \includegraphics[width=75mm]{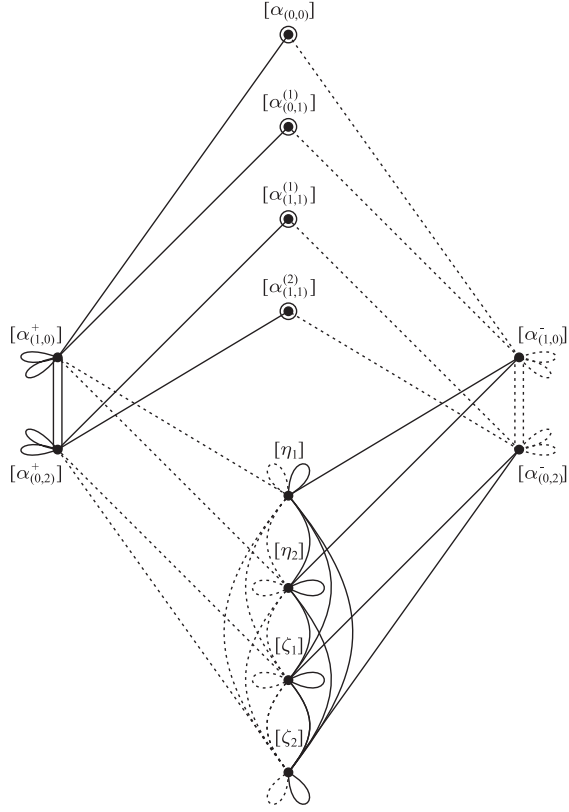} \\
 \caption{$\mathcal{E}_4(G_2)$: Multiplication by $[\alpha_{(1,0)}^+]$ (solid lines) and $[\alpha_{(1,0)}^-]$ (dashed lines)} \label{Fig-full_E4_G2}
\end{center}
\end{figure}

As in Section \ref{sect:measures_E3G2}, we can compute the principal graph and dual principal graph of the inclusion $(G_2)_4 \rightarrow (D_7)_1$. The chiral induced sector bases ${}_M \mathcal{X}_M^{\pm}$ and full induced sector basis ${}_M \mathcal{X}_M$ are given by
\begin{align*}
{}_M \mathcal{X}_M^{\pm} &= \{ [\alpha_{(0,0)}], [\alpha_{(1,0)}^{\pm}], [\alpha_{(0,2)}^{\pm}], [\alpha_{(0,1)}^{(1)}], [\alpha_{(1,1)}^{(1)}], [\alpha_{(1,1)}^{(2)}] \},\\
{}_M \mathcal{X}_M &= \{ [\alpha_{(0,0)}], [\alpha_{(1,0)}^+], [\alpha_{(1,0)}^-], [\alpha_{(0,2)}^+], [\alpha_{(0,2)}^-], [\alpha_{(0,1)}^{(1)}], [\alpha_{(1,1)}^{(1)}], [\alpha_{(1,1)}^{(2)}], [\eta_1], [\eta_2], [\zeta_1], [\zeta_2] \},
\end{align*}
where $[\alpha_{(0,1)}^{\pm}] = [\alpha_{(0,2)}^{\pm}] \oplus [\alpha_{(0,1)}^{(1)}]$, $[\alpha_{(1,1)}^{\pm}] = [\alpha_{(1,0)}^{\pm}] \oplus [\alpha_{(0,2)}^{\pm}] \oplus [\alpha_{(1,1)}^{(1)}] \oplus [\alpha_{(1,1)}^{(2)}]$, $[\alpha_{(1,0)}^+ \circ \alpha_{(1,0)}^-] = [\eta_1] \oplus [\eta_2]$ and $[\alpha_{(1,0)}^+ \circ \alpha_{(0,2)}^-] = [\zeta_1] \oplus [\zeta_2]$.
The fusion graphs of $[\alpha_{(1,0)}^+]$ (solid lines) and $[\alpha_{(1,0)}^-]$ (dashed lines) are given in Figure \ref{Fig-full_E4_G2}, where we have circled the marked vertices, and we note again that multiplication by $[\alpha_{(1,0)}^+]$ (or $[\alpha_{(1,0)}^-]$) gives one copy each of $\mathcal{E}_4(G_2)$ and $\mathcal{E}_4^M(G_2)$. The ambichiral part ${}_M \mathcal{X}_M^0$ obeys $\mathbb{Z}_4$ fusion rules, corresponding to $D_7$ at level 1.

\begin{figure}[tb]
\begin{center}
  \includegraphics[width=94mm]{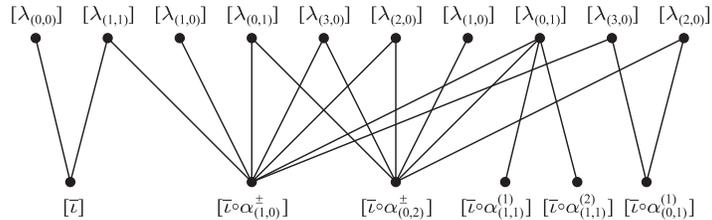} \\
 \caption{$\mathcal{E}_4(G_2)$: Principal graph of $(G_2)_4 \rightarrow (D_7)_1$} \label{Fig-GHJ_Graph_E4_G2}
\end{center}
\end{figure}

\begin{figure}[tb]
\begin{center}
  \includegraphics[width=92mm]{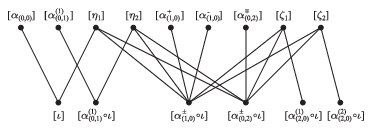} \\
 \caption{$\mathcal{E}_4(G_2)$: Dual principal graph of $(G_2)_4 \rightarrow (D_7)_1$} \label{Fig-GHJ_Graph_E4_G2_dual}
\end{center}
\end{figure}

We find
\begin{equation}
[\gamma] = [\alpha_{(0,0)}] \oplus [\eta_1],
\end{equation}
and the principal graph and dual principal graph of the inclusion $(G_2)_4 \rightarrow (D_7)_1$ of index $6+2\sqrt{6} \approx 10.90$ are illustrated in Figures \ref{Fig-GHJ_Graph_E4_G2} and \ref{Fig-GHJ_Graph_E4_G2_dual} respectively.

\begin{figure}[tb]
\begin{center}
  \includegraphics[width=45mm]{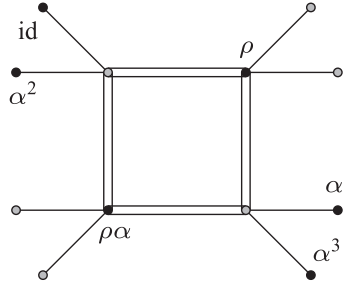} \\
 \caption{$\mathcal{E}_4(G_2)$: Principal graph of $\alpha_{(1,0)}^{\pm}(M) \subset M$} \label{Fig-principal_graph_E4_G2}
\end{center}
\end{figure}

Again, we can construct a subfactor $\alpha_{(1,0)}^{\pm}(M) \subset M$ where $M$ is a type III factor. Here the chiral systems ${}_M \mathcal{X}_M^{\pm}$ give quadratic extensions of a group category. Such fusion categories are discussed in \cite[$\S$1]{evans/gannon:2012}. In our case, $G$ is $\mathbb{Z}_4$ with subgroup $N = \mathbb{Z}_2$, $\rho = \alpha_{(1,0)}^{\pm}$ and $g_{\rho} = \mathrm{id}$ since $\rho$ is self-conjugate. The fusion rules are
\begin{eqnarray}
& [\alpha][\rho] = [\rho][\alpha] = [\rho \alpha] \neq [\rho], \qquad [\alpha^2][\rho] = [\rho] = [\rho][\alpha^2], & \label{eqn:fusion_rules-quad_group_cat-1} \\
& [\rho]^2 = 2[\rho] \oplus 2[\rho \alpha] \oplus [\mathrm{id}] \oplus [\alpha^2], & \label{eqn:fusion_rules-quad_group_cat-2}
\end{eqnarray}
where $[\alpha]$ satisfies $\mathbb{Z}_4$ fusion rules, $[\mathrm{id}]=[\alpha_{(0,0)}]$, $[\alpha]=[\alpha_{(1,1)}^{(1)}]$, $[\alpha^2]=[\alpha_{(0,1)}^{(1)}]$, $[\alpha^3]=[\alpha_{(1,1)}^{(2)}]$ and $[\rho \alpha] = [\alpha_{(0,2)}^{\pm}]$. Since $d_{\alpha}=1$, the index $d_{\rho}^2$ of $\rho(M) \subset M$ satisfies $d_{\rho}^2 = 4d_{\rho}+2$, thus $d_{\rho}^2 = (2+\sqrt{6})^2 = 10+4\sqrt{6} \approx 19.80$.
Its principal graph is illustrated in Figure \ref{Fig-principal_graph_E4_G2} and is the bipartite unfolded version of the graph $\mathcal{E}_4^{\rho_1}(G_2)$. The dual principal graph is again isomorphic to the principal graph as abstract graphs.

\begin{figure}[tb]
\begin{center}
  \includegraphics[width=55mm]{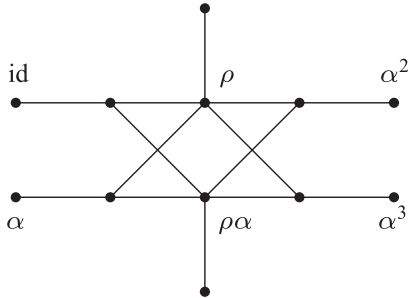} \\
 \caption{$\mathcal{E}_4(G_2)$: Principal graph of $\rho(M) \rtimes \mathbb{Z}_2 \subset M$} \label{Fig-quad_group_cat}
\end{center}
\end{figure}

For the fusion category above obtained from the conformal inclusion $(G_2)_4 \rightarrow (D_7)_1$, $\alpha_{(1,1)}^{(j)}$ is a non-trivial simple current of order 4 in ${}_M \mathcal{X}_M^0$, $j=1,2$. However, $[\alpha_{(1,1)}^{(j)}]$ is a subsector of $[\alpha_{(1,1)}^{\pm}]$, for which $\omega_{(1,1)} = e^{7\pi i/4}$ \cite{coquereaux/rais/tahri:2010}, thus $\omega_{(1,1)}^4 = e^{7\pi i} = -1$, and hence the orbifold inclusion $(G_2)_4 \rightarrow (D_7)_1 \rtimes_{\tau'} \mathbb{Z}_4$ does not exist (c.f. Section \ref{sect:measures_E3MG2}). On the other hand, the conformal dimension of the simple current $\alpha_{(0,1)}^{\pm}$ of order 2 is $1/2$ (mod $\mathbb{Z}$). Thus one can choose an automorphism $\beta$ on $M$ such that $[\beta] = [\alpha^2] \; (=[\alpha_{(0,1)}^{\pm}])$ and $\beta^2 = \mathrm{id}$. Thus there is an intermediate subfactor $\rho(M) \subset \rho(M) \rtimes \mathbb{Z}_2$ of index 2. The other intermediate subfactor $\rho(M) \rtimes \mathbb{Z}_2 \subset M$ would have index $5+2\sqrt{6}$. Writing the inclusions as $\iota: \rho(M) \rightarrow \rho(M) \rtimes \mathbb{Z}_2$ and $\jmath: \rho(M) \rtimes \mathbb{Z}_2 \rightarrow M$, as $M$-$M$ sectors the canonical endomorphism $\jmath \overline{\jmath}$ is a subsector of the canonical endomorphism $\jmath \iota \overline{\iota} \overline{\jmath} = \rho^2$. Hence $[\jmath \overline{\jmath}]$ is a subsector of $2[\rho] \oplus 2[\rho \alpha] \oplus [\mathrm{id}] \oplus [\beta]$ which contains $[\mathrm{id}]$. By considering the statistical dimensions, we see that $[\jmath \overline{\jmath}] = [\mathrm{id}] \oplus S$ for $S \in \{ 2[\rho], 2[\rho \alpha], [\rho] \oplus [\rho \alpha] \}$. The first two possibilities are not consistent with the fusion rules (\ref{eqn:fusion_rules-quad_group_cat-1})-(\ref{eqn:fusion_rules-quad_group_cat-2}), so we obtain $[\jmath \overline{\jmath}] = [\mathrm{id}] \oplus [\rho] \oplus [\rho \alpha]$, and the principal graph of $\rho(M) \rtimes \mathbb{Z}_2 \subset M$ is as in Figure \ref{Fig-quad_group_cat}.

We now determine the joint spectral measure of $\mathcal{E}_4^{\rho_1}(G_2)$, $\mathcal{E}_4^{\rho_2}(G_2)$.
With $\theta_1, \theta_2$ as in (\ref{eqn:theta-lambda_G2}) for $\lambda = (\lambda_1,\lambda_2) \in \mathrm{Exp}(\mathcal{E}_4(G_2))$, we have:
\begin{center}
\begin{tabular}{|c|c|c|c|} \hline
$\lambda \in \mathrm{Exp}$ & $(\theta_1,\theta_2) \in [0,1]^2$ & $|\psi^{\lambda}_{\ast}|^2$ & $\frac{1}{8\pi^2}|J(\theta_1,\theta_2)|$ \\
\hline $(0,0)$ & $\left(\frac{1}{6},\frac{23}{24}\right)$ & $\frac{3-\sqrt{6}}{24}$ & $\frac{3-\sqrt{6}}{\sqrt{3}}$ \\
\hline $(3,0)$ & $\left(\frac{7}{24},\frac{5}{6}\right)$ & $\frac{3+\sqrt{6}}{24}$ & $\frac{3+\sqrt{6}}{\sqrt{3}}$ \\
\hline $(0,1)$ & $\left(\frac{7}{24},\frac{23}{24}\right)$ & $\frac{1}{8}$ & $\sqrt{3}$ \\
\hline $(4,0)$ & $\left(\frac{1}{3},\frac{19}{24}\right)$ & $\frac{1}{8}$ & $\sqrt{3}$ \\
\hline $(1,1)$ & $\left(\frac{1}{3},\frac{11}{12}\right)$ & $0, \frac{1}{2}$ & $2\sqrt{3}$ \\
\hline
\end{tabular}
\end{center}
where again the eigenvectors $\psi^{\lambda}$ have been normalized so that $||\psi^{\lambda}|| = 1$. For the repeated exponent $(1,1)$, one of the eigenvectors has $|\psi^{(1,1)_1}_{\ast}|^2 = 0$ and the other has $|\psi^{(1,1)_2}_{\ast}|^2 = 1/2$, thus their sum $|\psi^{(1,1)_1}_{\ast}|^2 + |\psi^{(1,1)_2}_{\ast}|^2 = 1/2$.
Note that $24 |\psi^{\lambda}_{\ast}|^2 = \sqrt{3} |J|/8\pi^2$ for $\lambda \in \{ (0,0),(3,0) \}$, $3|\psi^{\lambda}_{\ast}|^2 = 2J^2/64\pi^4$ for $\lambda \in \{ (0,1),(4,0) \}$ and $24(|\psi^{(1,1)_1}_{\ast}|^2 + |\psi^{(1,1)_2}_{\ast}|^2) = J^2/64\pi^4$.

\begin{figure}[tb]
\begin{center}
  \includegraphics[width=45mm]{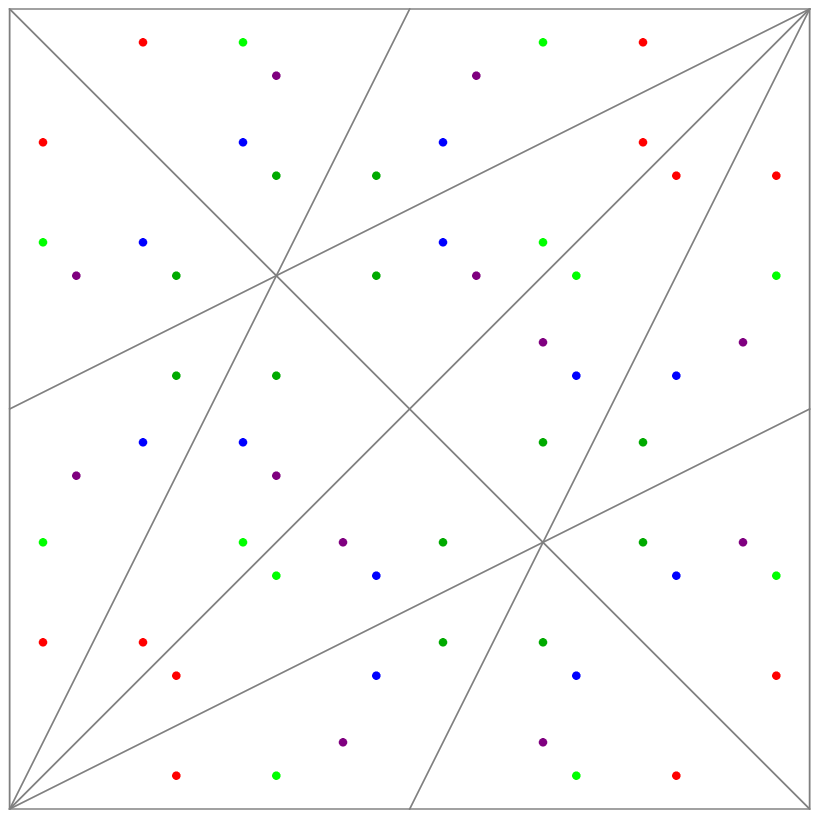}\\
 \caption{The orbit of the points $(\theta_1,\theta_2)$ for $\lambda \in \mathrm{Exp}(\mathcal{E}_4(G_2)$.} \label{Fig-E4pointsG2}
\end{center}
\end{figure}

The orbit under $D_{12}$ of $(\theta_1,\theta_2) \in \left\{ \left(\frac{1}{6},\frac{23}{24}\right)\textcolor{red}{\cdot}, \left(\frac{7}{24},\frac{5}{6}\right)\textcolor{blue}{\cdot}, \left(\frac{7}{24},\frac{23}{24}\right)\textcolor{green}{\cdot}, \left(\frac{7}{24},\frac{23}{24}\right)\textcolor{PineGreen}{\cdot}, \left(\frac{1}{3},\frac{11}{12}\right)\textcolor{purple}{\cdot} \right\}$ are illustrated in Figure \ref{Fig-E4pointsG2}.
Let $\Upsilon_{\lambda} := \sum_{g \in D_{12}} (\beta^{(g(\lambda))})^m (\beta^{(g(\lambda))})^n |\psi^{g(\lambda)}_{\ast}|^2$.
The orbits of the points $(1/6,23/24)$ and $(7/24,5/6)$ give the measure $\mathrm{d}^{(6,1/24)}$, whose support has cardinality 36, but we have to remove the additional points $\{ g(e^{\pi i/4},e^{\pi i}) | \, g \in D_{12} \}$ which are in the support of $\mathrm{d}^{(6,1/24)}$. For these additional points we have $|J| = 16\sqrt{2}\pi^2$.
Thus $\Upsilon_{(0,0)} + \Upsilon_{(3,0)} = \frac{36}{12} \frac{\sqrt{3}}{24} \frac{1}{8\pi^2}|J| \, \mathrm{d}^{(6,1/24)} - \frac{\sqrt{6}}{144} \sum_{g \in D_{12}} \delta_{g(e^{\pi i/4},-1)}$.
Now $\Upsilon_{(4,0)} = J^2/64\pi^4 \, \mathrm{d}^{((6))}$ since $J(\theta_1,\theta_2)=0$ for the additional points $(\theta_1,\theta_2) \in \mathrm{Supp}(\mathrm{d}^{((6))}) \setminus \{ g(1/6,19/24) | \, g \in D_{12} \} = \{ g(1/6,1/6) | \, g \in D_{12} \}$. We also have $\Upsilon_{(0,1)} = J^2/64\pi^4 \, \mathrm{d}^{((8/3))}$ since again $J(\theta_1,\theta_2)=0$ for the additional points in $\mathrm{Supp}(\mathrm{d}^{((8/3))})$, and similarly $\Upsilon_{(1,1)} = J^2/1024\pi^4 \, \mathrm{d}^{((4))}$.

Thus we have the following result:

\begin{Thm}
The joint spectral measure of $\mathcal{E}_4^{\rho_1}(G_2)$, $\mathcal{E}_4^{\rho_2}(G_2)$ (over $\mathbb{T}^2$) is
\begin{eqnarray}
\mathrm{d}\varepsilon & = & \frac{\sqrt{3}}{64\pi^2} |J| \, \mathrm{d}^{(6,1/24)} + \frac{1}{1024\pi^4} J^2 \, \mathrm{d}^{((4))} + \frac{1}{64\pi^4} J^2 \, \mathrm{d}^{((6))} + \frac{1}{64\pi^4} J^2 \, \mathrm{d}^{((8/3))} \nonumber \\
&& \quad - \frac{\sqrt{6}}{144} \sum_{g \in D_{12}} \delta_{g(e^{\pi i/4},-i)},
\end{eqnarray}
where $\mathrm{d}^{((n))}$, $\mathrm{d}^{(n,k)}$ are as in Definition \ref{def:4measures} and $\delta_x$ is the Dirac measure at the point $x$.
\end{Thm}

\subsection{Exceptional Graph $\mathcal{E}_4^M(G_2)$: $(G_2)_4 \rightarrow (D_7)_1 \rtimes \mathbb{Z}_2$}

\begin{figure}[tb]
\begin{minipage}[t]{5.5cm}
\begin{center}
  \includegraphics[width=35mm]{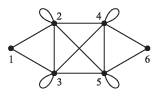} \\
 \caption{Graph $\mathcal{E}_4^{M,\rho_1}(G_2)$} \label{fig-Graph_E4M_G2-1}
\end{center}
\end{minipage}
\hfill
\begin{minipage}[t]{10cm}
\begin{center}
  \includegraphics[width=100mm]{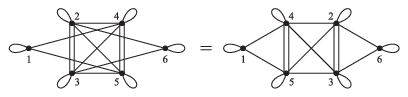} \\
 \caption{Graph $\mathcal{E}_4^{M,\rho_2}(G_2)$} \label{fig-Graph_E4M_G2-2}
\end{center}
\end{minipage}
\end{figure}

\begin{figure}[tb]
\begin{center}
  \includegraphics[width=90mm]{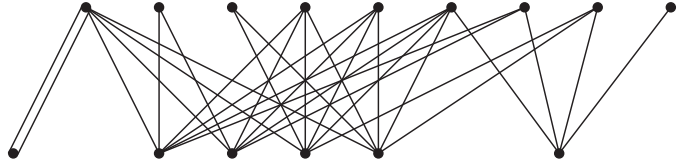} \\
 \caption{$\mathcal{E}_4^M(G_2)$: Principal graph of $(G_2)_4 \rightarrow (D_7)_1 \rtimes \mathbb{Z}_2$} \label{Fig-GHJ_Graph_E4M_G2}
\end{center}
\end{figure}

The graphs $\mathcal{E}_4^{M,\rho_j}(G_2)$, illustrated in Figures \ref{fig-Graph_E4M_G2-1} and \ref{fig-Graph_E4M_G2-2} are the nimrep graphs for the type II inclusion $(G_2)_4 \rightarrow (D_7)_1 \rtimes_{\tau} \mathbb{Z}_2$ with index $12+4\sqrt{6} \approx 21.80$, where $\tau = \alpha_{(0,1)}^{(1)}$ is a non-trivial simple current of order 2 in the ambichiral system ${}_M \mathcal{X}_M^0$, see Section \ref{sect:measures_E4G2}.
The principal graph for this inclusion is illustrated in Figure \ref{Fig-GHJ_Graph_E4M_G2}, which will be discussed in a future publication using a generalised Goodman-de la Harpe-Jones construction (c.f. the comments in Section \ref{sect:measures_E3MG2}). Again, it is not clear what the dual principal graph is in this case.

The associated modular invariant is again $Z_{\mathcal{E}_4}$ and the graphs are isospectral to $\mathcal{E}_4^{\rho_j}(G_2)$, with $\mathcal{E}_4^{M,\rho_j}(G_2)$ obtained from $\mathcal{E}_4^{\rho_j}(G_2)$ by a $\mathbb{Z}_2$-orbifold procedure.
However, the eigenvectors $\psi^{\lambda}$ are not identical to those for $\mathcal{E}_4^{\rho_j}(G_2)$.
With $\theta_1, \theta_2$ as in (\ref{eqn:theta-lambda_G2}) for $\lambda = (\lambda_1,\lambda_2) \in \mathrm{Exp}(\mathcal{E}_4^M(G_2)) = \mathrm{Exp}(\mathcal{E}_4(G_2))$, we have:
\begin{center}
\begin{tabular}{|c|c|c|c|} \hline
$\lambda \in \mathrm{Exp}$ & $(\theta_1,\theta_2) \in [0,1]^2$ & $|\psi^{\lambda}_{\ast}|^2$ & $\frac{1}{8\pi^2}|J(\theta_1,\theta_2)|$ \\
\hline $(0,0)$ & $\left(\frac{1}{6},\frac{23}{24}\right)$ & $\frac{3-\sqrt{6}}{60}$ & $\frac{3-\sqrt{6}}{\sqrt{3}}$ \\
\hline $(3,0)$ & $\left(\frac{7}{24},\frac{5}{6}\right)$ & $\frac{3+\sqrt{6}}{60}$ & $\frac{3+\sqrt{6}}{\sqrt{3}}$ \\
\hline $(0,1)$ & $\left(\frac{7}{24},\frac{23}{24}\right)$ & $\frac{1}{4}$ & $\sqrt{3}$ \\
\hline $(4,0)$ & $\left(\frac{1}{3},\frac{19}{24}\right)$ & $\frac{1}{4}$ & $\sqrt{3}$ \\
\hline $(1,1)$ & $\left(\frac{1}{3},\frac{11}{12}\right)$ & $0, 0$ & $2\sqrt{3}$ \\
\hline
\end{tabular}
\end{center}
where the eigenvectors $\psi^{\lambda}$ have been normalized so that $||\psi^{\lambda}|| = 1$.
In this case $60 |\psi^{\lambda}_{\ast}|^2 = \sqrt{3} |J|/8\pi^2$ for $\lambda \in \{ (0,0),(3,0) \}$ and $48|\psi^{\lambda}_{\ast}|^2 = 2J^2/64\pi^4$ for $\lambda \in \{ (0,1),(4,0) \}$.
Thus we have the following result:

\begin{Thm}
The joint spectral measure of $\mathcal{E}_4^{M,\rho_1}(G_2)$, $\mathcal{E}_4^{M,\rho_2}(G_2)$ (over $\mathbb{T}^2$) is
\begin{equation}
\mathrm{d}\varepsilon = \frac{\sqrt{3}}{160\pi^2} |J| \, \mathrm{d}^{(6,1/24)} + \frac{1}{2048\pi^4} J^2 \, \mathrm{d}^{((6))} + \frac{1}{2048\pi^4} J^2 \, \mathrm{d}^{((8/3))} - \frac{\sqrt{6}}{360} \sum_{g \in D_{12}} \delta_{g(e^{\pi i/4},-i)},
\end{equation}
where $\mathrm{d}^{((n))}$, $\mathrm{d}^{(n,k)}$ are as in Definition \ref{def:4measures} and $\delta_x$ is the Dirac measure at the point $x$.
\end{Thm}

\subsection{Exceptional Graph $\mathcal{E}_4^{\ast}(G_2)$}

\begin{figure}[tb]
\begin{minipage}[t]{7.5cm}
\begin{center}
  \includegraphics[width=35mm]{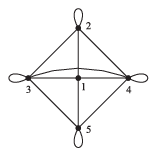} \\
 \caption{Graph $\mathcal{E}_4^{\ast,\rho_1}(G_2)$} \label{fig-Graph_E4star_G2-1}
\end{center}
\end{minipage}
\hfill
\begin{minipage}[t]{7.5cm}
\begin{center}
  \includegraphics[width=35mm]{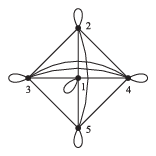} \\
 \caption{Graph $\mathcal{E}_4^{\ast,\rho_2}(G_2)$} \label{fig-Graph_E4star_G2-2}
\end{center}
\end{minipage}
\end{figure}

The graphs $\mathcal{E}_4^{\ast,\rho_j}(G_2)$ are illustrated in Figures \ref{fig-Graph_E4star_G2-1} and \ref{fig-Graph_E4star_G2-2}.
To our knowledge the second graph $\mathcal{E}_4^{\ast,\rho_2}(G_2)$ has not appeared in the literature before in the context of nimrep graphs or subfactors.
The associated modular invariant is \cite[(5.1)]{verstegen:1990}
\begin{align*}
Z_{\mathcal{E}_4^{\ast}} = |\chi_{(0,0)}|^2 &+ |\chi_{(3,0)}|^2 + |\chi_{(1,1)}|^2 + |\chi_{(2,0)}|^2 + |\chi_{(2,1)}|^2 \\
&+ \chi_{(1,0)}\chi_{(0,2)}^{\ast} + \chi_{(0,2)}\chi_{(1,0)}^{\ast} + \chi_{(0,1)}\chi_{(4,0)}^{\ast} + \chi_{(4,0)}\chi_{(0,1)}^{\ast}
\end{align*}
which is at level 4 and has exponents $\mathrm{Exp}(\mathcal{E}_4^{\ast}(G_2)) = \{ (0,0), (3,0), (1,1), (2,0), (2,1) \}$.

\begin{figure}[tb]
\begin{center}
  \includegraphics[width=105mm]{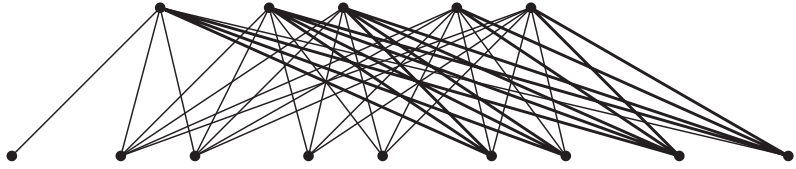} \\
 \caption{Expected principal graph of $G_2$-GHJ subfactor with nimrep $\mathcal{E}_4^{\ast}(G_2)$} \label{Fig-GHJ_Graph_E4star_G2}
\end{center}
\end{figure}

This modular invariant is a permutation invariant, and does not come from a conformal embedding. It has not yet been shown that the graphs $\mathcal{E}_4^{\ast,\rho_j}(G_2)$ arise from a braided subfactor. This will be discussed in a future publication using a generalised Goodman-de la Harpe-Jones construction (c.f. the comments in Section \ref{sect:measures_E3MG2}), which produces the second graph $\mathcal{E}_4^{\ast,\rho_2}(G_2)$ as a nimrep graph.
It is expected that $\mathcal{E}_4^{\ast,\rho_j}(G_2)$ does indeed arise as the nimrep for a type II inclusion with index $39+16\sqrt{6} \approx 78.19$. The expected principal graph for this inclusion is illustrated in Figure \ref{Fig-GHJ_Graph_E4star_G2}, where the thick lines indicate double edges. Again, it is not clear what the dual principal graph is in this case.

However, for our purposes it is sufficient to know the eigenvalues and corresponding eigenvectors for these graphs, and it is not necessary for the graph to be a nimrep graph.
For this graph there are two distinct vertices (up to an automorphism of the graph) which both have lowest Perron-Frobenius weight. These are numbered 1 and 2 in both Figures \ref{fig-Graph_E4star_G2-1}, \ref{fig-Graph_E4star_G2-2}.
Here we compute the spectral measure where the distinguished vertex is one of these vertices, the vertex numbered 1. Choosing the other vertex with lowest Perron-Frobenius weight as the distinguished vertex would yield a different measure.

Then with $\theta_1, \theta_2$ as in (\ref{eqn:theta-lambda_G2}) for $\lambda = (\lambda_1,\lambda_2) \in \mathrm{Exp}(\mathcal{E}_4^{\ast}(G_2))$, we have:
\begin{center}
\begin{tabular}{|c|c|c|} \hline
$\lambda \in \mathrm{Exp}$ & $(\theta_1,\theta_2) \in [0,1]^2$ & $|\psi^{\lambda}_{\ast}|^2$ \\
\hline $(0,0)$ & $\left(\frac{1}{6},\frac{23}{24}\right)$ & $\frac{1}{6}$ \\
\hline $(3,0)$ & $\left(\frac{7}{24},\frac{5}{6}\right)$ & $\frac{1}{6}$ \\
\hline $(1,1)$ & $\left(\frac{1}{3},\frac{11}{12}\right)$ & $0$ \\
\hline $(2,0)$ & $\left(\frac{1}{4},\frac{7}{8}\right)$ & $0$ \\
\hline $(2,1)$ & $\left(\frac{3}{8},\frac{7}{8}\right)$ & $\frac{2}{3}$ \\
\hline
\end{tabular}
\end{center}
where again the eigenvectors $\psi^{\lambda}$ have been normalized so that $||\psi^{\lambda}|| = 1$.
We have the following result:

\begin{Thm}
A spectral measure of $\mathcal{E}_4^{\ast,\rho_1}(G_2)$, $\mathcal{E}_4^{\ast,\rho_2}(G_2)$ (over $\mathbb{T}^2$), with distinguished vertex $\ast = 1$ in Figure \ref{fig-Graph_E4star_G2-1}, is
\begin{equation}
\mathrm{d}\varepsilon = \frac{1}{2} \mathrm{d}^{(6,1/24)} + \frac{1}{18} \sum_{g \in D_{12}} \delta_{g(i,e^{3\pi i/4})},
\end{equation}
where $\mathrm{d}^{(n,k)}$ is as in Definition \ref{def:4measures} and $\delta_x$ is the Dirac measure at the point $x$.
\end{Thm}

\bigskip

\begin{footnotesize}
\noindent{\it Acknowledgement.}

The second author was supported by the Coleg Cymraeg Cenedlaethol.
\end{footnotesize}

\end{document}